
%
\magnification1200
\pretolerance=100
\tolerance=200
\hbadness=1000
\vbadness=1000
\linepenalty=10
\hyphenpenalty=50
\exhyphenpenalty=50
\binoppenalty=700
\relpenalty=500
\clubpenalty=5000
\widowpenalty=5000
\displaywidowpenalty=50
\brokenpenalty=100
\predisplaypenalty=7000
\postdisplaypenalty=0
\interlinepenalty=10
\doublehyphendemerits=10000
\finalhyphendemerits=10000
\adjdemerits=160000
\uchyph=1
\delimiterfactor=901
\hfuzz=0.1pt
\vfuzz=0.1pt
\overfullrule=5pt
\hsize=146 true mm
\vsize=8.9 true in
\maxdepth=4pt
\delimitershortfall=.5pt
\nulldelimiterspace=1.2pt
\scriptspace=.5pt
\normallineskiplimit=.5pt
\mathsurround=0pt
\parindent=20pt
\catcode`\_=11
\catcode`\_=8
\normalbaselineskip=12pt
\normallineskip=1pt plus .5 pt minus .5 pt
\parskip=6pt plus 3pt minus 3pt
\abovedisplayskip = 12pt plus 5pt minus 5pt
\abovedisplayshortskip = 1pt plus 4pt
\belowdisplayskip = 12pt plus 5pt minus 5pt
\belowdisplayshortskip = 7pt plus 5pt
\normalbaselines
\smallskipamount=\parskip
 \medskipamount=2\parskip
 \bigskipamount=3\parskip
\jot=3pt
%
%
\def\ref#1{\par\noindent\hangindent2\parindent
 \hbox to 2\parindent{#1\hfil}\ignorespaces}
%
%
\font\typd=cmbx10 scaled \magstep2   
\font\typf=cmcsc10                   
\font\tenss=cmss10
\font\sevenss=cmss8 at 7pt
\font\fivess=cmss8 at 5pt
\newfam\ssfam %
\textfont\ssfam=\tenss
\scriptfont\ssfam=\sevenss
\scriptscriptfont\ssfam=\fivess
%
%
%
%
%
%
%
%
%
\catcode`\_=11
\def\suf_fix{}
\def\scaled_rm_box#1{%
 \relax
 \ifmmode
   \mathchoice
    {\hbox{\tenrm #1}}%
    {\hbox{\tenrm #1}}%
    {\hbox{\sevenrm #1}}%
    {\hbox{\fiverm #1}}%
 \else
  \hbox{\tenrm #1}%
 \fi}
\def\suf_fix_def#1#2{\expandafter\def\csname#1\suf_fix\endcsname{#2}}
\def\I_Buchstabe#1#2#3{%
 \suf_fix_def{#1}{\scaled_rm_box{I\hskip-0.#2#3em #1}}
}
\def\rule_Buchstabe#1#2#3#4{%
 \suf_fix_def{#1}{%
  \scaled_rm_box{%
   \hbox{%
    #1%
    \hskip-0.#2em%
    \lower-0.#3ex\hbox{\vrule height1.#4ex width0.07em }%
   }%
   \hskip0.50em%
  }%
 }%
}
\I_Buchstabe B22
\rule_Buchstabe C51{34}
\I_Buchstabe D22
\I_Buchstabe E22
\I_Buchstabe F22
\rule_Buchstabe G{525}{081}4
\I_Buchstabe H22
\I_Buchstabe I20
\I_Buchstabe K22
\I_Buchstabe L20
\I_Buchstabe M{20em }{I\hskip-0.35}
\I_Buchstabe N{20em }{I\hskip-0.35}
\rule_Buchstabe O{525}{095}{45}
\I_Buchstabe P20
\rule_Buchstabe Q{525}{097}{47}
\I_Buchstabe R21 
\rule_Buchstabe U{45}{02}{54}
\suf_fix_def{Z}{\scaled_rm_box{Z\hskip-0.38em Z}}
\catcode`\"=12
\newcount\math_char_code
\def\suf_fix_math_chars_def#1{%
 \ifcat#1A
  \expandafter\math_char_code\expandafter=\suf_fix_fam
  \multiply\math_char_code by 256
  \advance\math_char_code by `#1
  \expandafter\mathchardef\csname#1\suf_fix\endcsname=\math_char_code
  \let\next=\suf_fix_math_chars_def
 \else
  \let\next=\relax
 \fi
 \next}
%
%
%
%
\def\font_fam_suf_fix#1#2 #3 {%
 \def\suf_fix{#2}
 \def\suf_fix_fam{#1}
 \suf_fix_math_chars_def #3.
}
\font_fam_suf_fix
 0rm
 ABCDEFGHIJKLMNOPQRSTUVWXYZabcdefghijklmnopqrstuvwxyz
\font_fam_suf_fix
 2scr
 ABCDEFGHIJKLMNOPQRSTUVWXYZ
\font_fam_suf_fix
 \slfam sl
 ABCDEFGHIJKLMNOPQRSTUVWXYZabcdefghijklmnopqrstuvwxyz
\font_fam_suf_fix
 \bffam bf
 ABCDEFGHIJKLMNOPQRSTUVWXYZabcdefghijklmnopqrstuvwxyz
\font_fam_suf_fix
 \ttfam tt
 ABCDEFGHIJKLMNOPQRSTUVWXYZabcdefghijklmnopqrstuvwxyz
\font_fam_suf_fix
 \ssfam
 ss
 ABCDEFGHIJKLMNOPQRSTUVWXYZabcdefgijklmnopqrstuwxyz
\catcode`\_=8
\def\Adss{{\fam\ssfam A\mkern -9.5mu A}}%
\def\Cdss{{\fam\ssfam
    \mkern 4.2 mu \mathchoice%
    {\vrule height 6.5pt depth -.55pt width 1pt}%
    {\vrule height 6.5pt depth -.57pt width 1pt}%
    {\vrule height 4.55pt depth -.28pt width .8pt}%
    {\vrule height 3.25pt depth -.19pt width .6pt}%
    \mkern -6.3mu C}}%
\def\Gdss{{\fam\ssfam
    \mkern 3.8 mu \mathchoice%
    {\vrule height 6.5pt depth -.62pt width 1pt}%
    {\vrule height 6.5pt depth -.65pt width 1pt}%
    {\vrule height 4.55pt depth -.44pt width .8pt}%
    {\vrule height 3.25pt depth -.30pt width .6pt}%
    \mkern -5.9mu G}}%
\def\Ndss{{\fam\ssfam I\mkern -2.5mu N}}%
\def\Qdss{{\fam\ssfam
    \mkern 3.8 mu \mathchoice%
    {\vrule height 6.5pt depth -.67pt width 1pt}%
    {\vrule height 6.5pt depth -.7pt width 1pt}%
    {\vrule height 4.55pt depth -.44pt width .7pt}%
    {\vrule height 3.25pt depth -.3pt width .5pt}%
    \mkern -5.9mu Q}}%
\def\Zdss{{\fam\ssfam Z\mkern-8.1mu Z}}%
%
%
%
%
\font\teneuf=eufm10 
\font\seveneuf=eufm7
\font\fiveeuf=eufm5
\newfam\euffam \def\euf{\fam\euffam\teneuf} 
\textfont\euffam=\teneuf \scriptfont\euffam=\seveneuf
\scriptscriptfont\euffam=\fiveeuf
       \def\afr{{\euf a}}
       \def\bfr{{\euf b}}
       \def\cfr{{\euf c}}

       \def\ffr{{\euf f}}
       \def\gfr{{\euf g}}

       \def\mfr{{\euf m}}

       \def\xfr{{\euf x}}

\input xypic.tex
\parindent=0pt
\let\ul=\underline
\def\dQ{{\Qdss_p}}
\def\dC{{\Cdss_p}}
\def\dZ{{\Zdss_p}}
\def\GL{{\rm GL}}
\def\Hom{\mathop{\rm Hom}\nolimits}
\def\dlongrightarrow{\longrightarrow\hskip-8pt\rightarrow}
\def\llongmapsto{\longleftarrow\hskip-4.2pt\hbox{\vrule height5pt width 0.03em}\ }
\def\bfB{{\bf B}}
\def\dG{{\Gdss}_{m}}
\def\thfill{\null\nobreak\hfill}
\def\endbeweis{\thfill\vbox{\hrule
  \hbox{\vrule\hbox to 5pt{\vbox to 5pt{\vfil}\hfil}\vrule}\hrule}}
\def\frac#1#2{{{#1}\over{#2}}}
\def\into{{\hookrightarrow}}
\def\ord{\mathop{\rm ord}\nolimits}
\def\ophi{\overline{\varphi}}
\def\bfK{{\bf K}}
\def\bfF{{\bf F}}

\centerline{\typd p-adic Fourier theory}

\medskip

\centerline{\typf P. Schneider, J. Teitelbaum}

\bigskip

In the early sixties, Amice ([Am1], [Am2]) studied the space of
$K$-valued, locally analytic functions on $\dZ$ and formulated a
complete description of its dual, the ring of $K$-valued, locally
$\dQ$-analytic distributions on $\dZ$, when $K$ is a complete
subfield of $\dC$. She found an isomorphism between the ring of
distributions and the space of global functions on a rigid variety
over $K$  parameterizing $K$-valued, locally analytic characters
of $\dZ$.  This rigid variety is in fact the open unit disk,  a
point $z$ of $\dC$ with $|z|<1$ corresponding to the locally
$\dQ$-analytic character $\kappa_{z}(a)=(1+z)^a$ for $a\in\dZ$.
The rigid function $F_{\lambda}$ corresponding to a distribution
$\lambda$ is determined by the formula
$F_{\lambda}(z)=\lambda(\kappa_{z})$.  Amice's description of the
ring of $\dQ$-analytic distributions was complemented by results
of Lazard ([Laz]). He described a divisor theory for the ring of
functions on the open disk and proved that, when $K$ is
spherically complete, the classes of closed, finitely generated,
and principal ideals in this ring coincide.

In this paper we generalize the work of Amice and Lazard by
studying the space $C^{an}(o,K)$ of $K$-valued, locally
$L$-analytic functions on $o$, and the corresponding ring of
distributions $D(o,K)$, when  $\dQ \subseteq L \subseteq K
\subseteq \dC$ with $L$ finite over $\dQ$ and $K$ complete and $o = o_L$
the additive group of the ring of integers in $L$. To clarify
this, observe that, as a $\dQ$-analytic manifold, the ring $o$ is
a product of $[L:\dQ]$ copies of $\dZ$. The $K$-valued,
$\dQ$-analytic functions on $o$ are thus given locally by power
series in $[L:\dQ]$ variables, with coefficients in $K$. The
$L$-analytic functions in $C^{an}(o,K)$ are given locally by power
series in {\it one} variable; they form a subspace of the
$\dQ$-analytic functions cut out by a set of ``Cauchy-Riemann''
differential equations. These facts are treated in Section 1.

Like Amice, we develop a Fourier theory for the locally
$L$-analytic functions on $o$.  We construct (Section 2) a rigid
group variety $\widehat{o}$, defined over $L$, whose closed points
$z$ in a field $K$ parameterize $K$-valued locally $L$-analytic
characters $\kappa_{z}$ of $o$.  We then show that, for $K$ a
complete subfield of $\dC$, the ring of rigid functions on
$\widehat{o}/K$ is isomorphic to the ring $D(o,K)$, where the
isomorphism $\lambda\mapsto F_{\lambda}$ is defined by
$\lambda(\kappa_{z})=F_{\lambda}(z)$, just as in Amice's
situation.

The most novel aspect of this situation is the variety
$\widehat{o}$. We prove (Section 3) that $\widehat{o}$ is a rigid
variety defined over $L$ that becomes isomorphic over $\dC$ to the
open unit disk, but is not isomorphic to a disk over any
discretely valued extension field of $L$. The ring of rigid
functions on $\widehat{o}$ has the property that the classes of
closed, finitely generated, and invertible ideals coincide; but we
show that unless $L=\dQ$ (Lazard's situation) there are
non-principal, finitely generated ideals, even over spherically
complete coefficient fields.

The ``uniformization'' of  $\widehat{o}$ by the open unit disk
follows from a result of Tate's in his famous paper on
$p$-divisible groups ([Tat]).  We show that over $\dC$ the  group
$\widehat{o}$ becomes isomorphic to the group of $\dC$-valued
points of a Lubin-Tate formal group associated to $L$.   The
Galois cocycle that gives the descent data on the open unit disk
yielding the twisted form $\widehat{o}$ comes directly out of the
Lubin-Tate group. The period of the Lubin-Tate group  plays a
crucial role in the explicit form of our results; in an Appendix
we use results of Fontaine [Fon] to obtain information on the
valuation of this period, generalizing work of Boxall ([Box]).

We give two applications of our Fourier theory.  The first is a
generalized Mahler expansion for locally $L$-analytic functions on
$o$ (Section 4).  The second is a construction of a $p$-adic
L-function for a CM elliptic curve at a supersingular prime
(Section 5). Although the method by which we obtain it is more
natural, and we obtain stronger analyticity results,  the
L-function we construct is essentially that studied by Katz
([Kat]) and by Boxall ([Box]).  The paper by Katz in particular
was a major source of inspiration in our work.

Our original motivation for studying this problem came from our
work on locally analytic representation theory.  In the paper [ST]
we classified the locally analytic principal series
representations of $\GL_{2}(\dQ)$.  The results of Amice and
Lazard played a key role in the proof, and in seeking to
generalize those results to the principal series of $\GL_{2}(L)$
we were led to consider the problems discussed in this paper.  The
results of this paper are sufficient to extend the methods of [ST]
to the groups $\GL_{2}(L)$, though to keep the paper
self-contained we do not give the proof here.

The relationship between formal groups and $p$-adic integration
has been known and exploited in some form by many authors.  We
have already mentioned the work of Katz [Kat] and Boxall [Box].
Height one formal groups and their connection to $p$-adic
integration is systematically used in [dS] and we have adapted
this approach to the height two case in Section 5 of our paper.
Some other results of a similar flavor were obtained in [SI].
Finally, we  point out that the appearance of $p$-adic Hodge
theory in our work raises the interesting question of relating our
results to the work of Colmez ([Col]).

We would like to thank John Coates, Robert Coleman, Pierre Colmez,
and Ehud deShalit for interesting discussions on this work. The
first, resp. second, author thanks the Landau Center, resp. the
Forcheimer Fellowship Fund, and Hebrew University for its support
and hospitality during the writing of this paper. The second
author was also supported by a grant from the NSA Mathematical
Sciences Program.

\bigskip

{\bf 1. Preliminaries on restriction of scalars}

\smallskip

We fix fields $\dQ \subseteq L \subseteq K$ such that $L/\dQ$ is
finite and $K$ is complete with respect to a nonarchimedean
absolute value $|\ |$ extending the one on $L$. We also fix a
commutative $d$-dimensional locally $L$-analytic group $G$. Then
the locally convex $K$-vector space $C^{an} (G,K)$ of all
$K$-valued locally analytic functions on $G$ is defined ([Fe2]
2.1.10).

We consider now an intermediate field $\dQ \subseteq L_0 \subseteq
L$ and let $G_0$ denote the locally $L_0$-analytic group obtained
from $G$ by restriction of scalars ([B-VAR] \S5.14). The dimension
of $G_0$ is $d[L:L_0]$. There is the obvious injective continuous
$K$-linear map
$$
C^{an} (G,K) \longrightarrow C^{an} (G_0 ,K)\ .\leqno{(\ast)}
$$
We want to describe the image of this map. The Lie algebra $\gfr$
of $G$ can naturally be identified with the Lie algebra of $G_0$
([B-VAR] 5.14.5). We fix an exponential map $\exp:\gfr\, {--->}\,
G$ for $G$; it, in particular, is a local isomorphism, and can be
viewed as an exponential map for $G_0$ as well. The Lie algebra
$\gfr$ acts in a compatible way on both sides of the above map by
continuous endomorphisms defined by
$$
  (\xfr f)(g):={d\over dt} f(\exp (-t\xfr ) g)_{|t=0}
$$
([Fe2] 3.1.2 and 3.3.4). By construction the map $\xfr \rightarrow
\xfr f$ on $\gfr$, for a fixed $f \in C^{an} (G_0 ,K)$, resp. $f \in C^{an}
(G,K)$, is $L_0$-linear, resp. $L$-linear.

\medskip

{\bf Lemma 1.1:}

{\it The image of $(\ast)$ is the closed subspace of all $f \in
C^{an} (G_0 ,K)$ such that $(t\xfr )f = t\cdot (\xfr f)$ for any
$\xfr
\in \gfr$ and any $t \in L$.}

Proof: We fix an $L$-basis $\xfr_1 ,\ldots,\xfr_d$ of $\gfr$ as
well an orthonormal basis $v_1 = 1,v_2,\ldots,v_e$ of $L$ as a
normed $L_0$-vector space. Then $v_1 \xfr_1 ,\ldots,v_e \xfr_d$ is
an $L_0$-basis of $\gfr$. Using the corresponding canonical
coordinates of the second kind ([B-GAL] Chap.III,\S4.3) we have,
for a given $f \in C^{an} (G_0 ,K)$ and a given $g \in G$, the
convergent expansion
$$
f({\rm exp}(t_{11}v_1 \xfr_1 + \ldots + t_{ed}v_e \xfr_d)g) =
\sum_{n_{11},\ldots,n_{ed} \geq 0} c_{\ul n}t_{11}^{n_{11}}\ldots
t_{ed}^{n_{ed}}\ ,
$$
with $c_{\ul n} \in K$, in a neighborhood of $g$ (i.e., for
$t_{ij}\in L_0$ small enough). We now assume that
$$
(v_i \xfr_j)f = v_i \cdot (\xfr_j f) = v_i \cdot ((v_1 \xfr_j)f)
$$
holds true for all $i$ and $j$. Computing both sides in terms of
the above expansion and comparing coefficients results in the
equations
$$
(n_{ij}+1)c_{(n_{11},\ldots,n_{ij}+1,\ldots,n_{ed})} = v_i
(n_{1j}+1)c_{(n_{11},\ldots,n_{1j}+1,\ldots,n_{ed})}\ .
$$
Introducing the tuple ${\ul m}({\ul n}) = (m_1 ,\ldots,m_d)$
defined by $m_j := n_{1j} + \ldots + n_{ej}$ and the new
coefficients $b_{{\ul m}({\ul n})} := c_{(m_1 ,0,\ldots,m_2
,0,\ldots,m_d ,0,\ldots)}$ we deduce from this by induction that
$$
c_{\ul n} = b_{{\ul m}({\ul n})}{m_1 ! \over n_{11}!\ldots
n_{e1}!} \ldots {m_d ! \over n_{1d}!\ldots
n_{ed}!}v_1^{n_{11}+\ldots n_{1d}}\ldots v_e^{n_{e1}+\ldots
+n_{ed}}\ .
$$
Inserting this back into the above expansion and setting $t_j :=
t_{1j}v_1 + \ldots + t_{ej}v_e$ we obtain the new expansion
$$
f({\rm exp}(t_1 \xfr_1 + \ldots + t_d \xfr_d)g) = \sum_{m_1
,\ldots,m_d \geq 0} b_{\ul m}t_1^{m_1}\ldots t_d^{m_d}
$$
which shows that $f$ is locally analytic on $G$.

\medskip

{\bf Lemma 1.2:}

{\it The map $(\ast)$ is a homeomorphism onto its (closed) image.}

Proof: Let $H \subseteq G$ be a compact open subgroup. According
to [Fe2] 2.2.4) we then have
$$
C^{an} (G,K) =\prod\limits_{g\in G/H} C^{an} (H ,K) \ .
$$
A corresponding decomposition holds for $G_0$. This shows that it
suffices to consider the case where $G$ is compact. In this case
$(\ast)$ is a compact inductive limit of isometries between Banach
spaces ([Fe2] 2.3.2), and the assertion follows from [GKPS]
3.1.16.

\medskip

The continuous dual $D(G,K) := C^{an} (G,K)'$ is the algebra of
$K$-valued distributions on $G$. The multiplication is the
convolution product $\ast$ ([Fe1] 4.4.2 and 4.4.4).

We assume from now on that $G$ is compact. To describe the correct
topology on $D(G,K)$ we need to briefly recall the construction of
$C^{an}(G,K)$. Let $G \supseteq H_0 \supseteq H_1
\supseteq\ldots\supseteq H_n \supseteq\ldots$ be a fundamental
system of open subgroups such that each $H_n$ corresponds under
the exponential map to an $L$-affinoid disk. We then have, for
each $g\in G$ and $n\in\Ndss$, the $K$-Banach space
$\Fscr_{gH_n}(K)$ of $K$-valued $L$-analytic functions on the
coset $gH_n$ viewed as an $L$-affinoid disk. The space
$C^{an}(G,K)$ is the locally convex inductive limit
$$
C^{an}(G,K) = \lim\limits_{{\longrightarrow}\atop{n}}\,
\Fscr_n(G,K)
$$
of the Banach spaces
$$
\Fscr_n(G,K) := \prod_{g\in G/H_n} \Fscr_{gH_n}(K)\ .
$$
The dual $D(G,K)$ therefore coincides as a vector space with the
projective limit
$$
D(G,K) = \lim\limits_{{\longleftarrow}\atop{n}}\,
\Fscr_{n}(G,K)'
$$
of the dual Banach spaces. We always equip $D(G,K)$ with the
corresponding projective limit topology. (Using [GKPS] 3.1.7(vii)
and the open mapping theorem one can show that this topology in
fact coincides with the strong dual topology.) In particular,
$D(G,K)$ is a commutative $K$-Fr\'echet algebra. The dual of the
map $(\ast)$ is a continuous homomorphism of Fr\'echet algebras
$$
D(G_0 ,K) \dlongrightarrow D(G,K)\ .\leqno(\ast)'
$$
It is surjective since $C^{an}(G_0,K)$ as a compact inductive
limit is of countable type ([GKPS] 3.1.7(viii)) and hence
satisfies the Hahn-Banach theorem ([Sh] 4.2 and 4.4). By the open
mapping theorem $(\ast)'$ then is a quotient map.

The action of $\gfr$ on $C^{an} (G_0 ,K)$ induces an action of
$\gfr$ on $D(G_0 ,K)$ by $(\xfr \lambda)(f) := \lambda(-\xfr f)$.
This action is related to the algebra structure through the
$L_0$-linear inclusion
$$
\matrix{\iota : & \gfr & \longrightarrow & D(G_0 ,K)\hfill\cr
& \xfr & \longmapsto & [f \mapsto (-\xfr (f))(1)]}
$$
which satisfies
$$
\iota(\xfr)\ast\lambda = \xfr \lambda\ \ \hbox{for}\ \xfr \in \gfr\
\hbox{and}\ \lambda \in D(G_0 ,K)
$$
(see the end of section 2 in [ST]). Followed by $(\ast)'$ this
inclusion becomes $L$-linear.







Let $\widehat{G_0}(K) \subseteq C^{an} (G_0 ,K)$ denote the subset
of all $K$-valued locally analytic characters on $G_0$. Any $\chi
\in \widehat{G_0}(K)$ induces the $L_0$-linear map
$$
\matrix{d\chi : & \gfr & \longrightarrow & K\hfill\cr\cr
& \xfr & \longmapsto & {d \over dt}\chi (\exp(t\xfr))_{|t=0}\ .}
$$

\medskip

{\bf Lemma 1.3:}

$\widehat{G}(K) = \{\chi \in \widehat{G_0}(K) : d\chi\ \hbox{is}\
L\hbox{-linear}\}.$

Proof: Because of
$$
(-\xfr\chi)(g) = \chi (g)\cdot d\chi (\xfr)
$$
this is a consequence of Lemma 1.1.

\medskip

The lemma says that the diagram
$$
\xymatrix{
  \widehat{G}(K) \ar[d]_{d} \ar[r]^{\subseteq}
                & \widehat{G_0}(K) \ar[d]^{d}  \\
  \Hom_L (\gfr ,K)  \ar[r]^{\subseteq}
                & \Hom_{L_0}(\gfr ,K)}
$$
is cartesian.

We suppose from now on that $K$ is a subfield of $\dC$ ($=$ the
completion of an algebraic closure of $\dQ$). There is the natural
strict inclusion
$$
\matrix{
D(G_0 ,K) & = & \lim\limits_{{\longleftarrow}\atop{n}}\,
\Hom_K^{cont}(\Fscr_n(G_0,K),K)\hfill\cr & & \big\downarrow\cr\cr
D(G_0 ,\dC) & = & \lim\limits_{{\longleftarrow}\atop{n}}\,
\Hom_{\dC}^{cont}(\Fscr_n(G_0,\dC),\dC)\ . }
$$

The ${\ul{Fourier\ transform}}$ $F_{\lambda}$ of a $\lambda \in
D(G_0 ,K)$, by definition, is the function
$$
\matrix{F_{\lambda} : & \widehat{G_0}(\dC) & \longrightarrow &
\dC\hfill\cr & \chi & \longmapsto & \lambda (\chi)\ .}
$$

\medskip

{\bf Proposition 1.4:}

{\it i. For any $\lambda \in D(G_0 ,K)$ we have $\lambda = 0$ if
and only if $F_{\lambda} = 0$;

ii. $F_{\mu \ast \lambda} = F_{\mu}F_{\lambda}$ for any two $\mu
,\lambda \in D(G_0 ,K)$.}

Proof: [Fe1] Thm. 5.4.8 (recall that $G$ is assumed to be
compact). For the convenience of the reader we sketch the proof of
the first assertion in the case of the additive group $G_0 = G =
o_L$. (This is the only case in which we actually will use this
result in the next section. Moreover, the general proof is just an
elaboration of this special case.) Let $\bfr \subseteq o_L$ be an
arbitrary nonzero ideal viewed as an additive subgroup. We use the
convention to denote by $f|a+\bfr$, for any function $f$ on $o_L$
and any coset $a+\bfr \subseteq o_L$, the function on $o_L$ which
is equal to $f$ on the coset $a+\bfr$ and which vanishes
elsewhere. Suppose now that $F_{\lambda} = 0$, i.e., that
$\lambda(\chi) = 0$ for any $\chi\in \widehat{G}(\dC)$. Using the
character theory of finite abelian groups one easily concludes
that
$$
\lambda(\chi|a+\bfr) = 0\ \ \ \ \hbox{for any}\
\chi\in\widehat{G}(\dC)\ \hbox{and any coset}\ a+\bfr \subseteq
o_L\ .
$$
We apply this to the character $\chi_y(x) := {\rm exp}(yx)$ where
$y\in o_{\dC}$ is small enough and obtain by continuity that
$$
0 = \lambda(\chi_y|a+\bfr) = \sum_{n \geq 0} {y^n \over
n!}\lambda(x^n|a+\bfr)\ .
$$
Viewing the right hand side as a power series in $y$ in a small
neighbourhood of zero it follows that
$$
\lambda(x^n|a+\bfr) = 0\ \ \ \ \hbox{for any}\ n \geq 0\ \hbox{and
any coset}\ a+\bfr \subseteq o_L\ .
$$
Again from continuity we see that $\lambda = 0$.

\medskip

{\bf Corollary 1.5:}

{\it i. $\widehat{G}(\dC) = \{\chi \in \widehat{G_0}(\dC) :
F_{\iota(t\xfr)-t\iota (\xfr)}(\chi) = 0\ \hbox{\rm{for\ any}}\
\xfr
\in \gfr\ \hbox{\rm{and}}\ t \in L\}$;

ii. the kernel of $(\ast)'$ is the ideal $I(G) := \{\lambda \in
D(G_0 ,K) : F_{\lambda}|\widehat{G}(\dC) = 0\}$.}

Proof: The assertion i. is a consequence of Lemma 1.3 and the
identity
$$
F_{\iota(t\xfr)-t\iota (\xfr)}(\chi) = ((-t\xfr)\chi)(1) - t(-
\xfr\chi)(1) = d\chi (t\xfr) - t\cdot d\chi (\xfr)\ .
$$
The assertion ii. follows from Prop. 1.4.i (applied to $G$).


\bigskip

{\bf 2. The Fourier transform for $G = o_L$}

\smallskip

Let $\dQ \subseteq L \subseteq K \subseteq \dC$ again be a chain of
complete fields and let $o := o_L$ denote the ring of integers in $L$.
The aim of this section is to determine the image of the Fourier
transform for the compact additive group $G := o$. The restriction of
scalars $G_0$ will always be understood with respect to the extension
$L/\dQ$.

First we have to discuss briefly a certain way to write rigid
analytic polydisks in a coordinate free manner. Let $\bfB_1$
denote the rigid $L$-analytic open disk of radius one around the
point $1 \in L$; its $K$-points are $\bfB_1 (K) = \{z \in K :
|z-1|<1\}$. We note that the group $\dZ$ acts on $\bfB_1$ via the
rigid analytic automorphisms
$$
\matrix{\dZ\ \times\ \bfB_1 & \longrightarrow & \bfB_1\hfill\cr\cr
(a,z) & \longmapsto & z^{a} := \mathop{\sum}\limits_{n\geq 0}
{a\choose n}(z-1)^n}
$$
(compare [Sch] \S\S32 and 47). Hence, given any free $\dZ$-module
$M$ of finite rank $r$, we can in an obvious sense form the rigid
$L$-analytic variety $\bfB_1 \otimes_{\dZ} M$ whose $K$-points are
$\bfB_1 (K) \otimes_{\dZ} M$. Any choice of an $\dZ$-basis of $M$
defines an isomorphism between $\bfB_1
\otimes_{\dZ} M$ and an $r$-dimensional open
polydisk over $L$. In particular, the family of all affinoid
subdomains in $\bfB_1 \otimes_{\dZ} M$ has a countable cofinal
subfamily. Writing the ring $\Oscr (\bfB_1 \otimes_{\dZ} M)$ of
global holomorphic functions on $\bfB_1 \otimes_{\dZ} M$ as the
projective limit of the corresponding affinoid algebras we see
that $\Oscr (\bfB_1 \otimes_{\dZ} M)$ in a natural way is an
$L$-Fr\'echet algebra.

After this preliminary discussion we recall that the maps
$$
\matrix{\widehat{\dZ}(K) & \longleftrightarrow & \bfB_1 (K)\cr\cr
\chi & \longmapsto & \chi(1)\cr
\chi_z (a) := z^a & \llongmapsto & z}
$$
are bijections inverse to each other (compare [Am2] 1.1 and [B-GAL]
III.8.1). They straightforwardly generalize to the bijection
$$
\matrix{\bfB_1 (K) \otimes_{\dZ} \Hom_{\dZ}(o,\dZ) &
\mathop{\longrightarrow}\limits^{\sim} & \widehat{G_0}(K)\hfill\cr\cr
z \otimes \beta & \longmapsto & \chi_{z\otimes \beta}(g) :=
z^{\beta(g)}\ .}
$$
By transport of structure the right hand side therefore can and will
be considered as the $K$-points of a rigid analytic group variety
$\widehat{G_0}$ over $L$ (which is non-canonically isomorphic to an
open polydisk of dimension $[L:\dQ]$). By construction the Lie algebra
of $\widehat{G_0}$ is equal to $\Hom_{\dQ}(\gfr,L)$. One easily checks
that
$$
d\chi_{z\otimes \beta} = \log(z)\cdot\beta\ .
$$
If we combine this identity with the commutative diagram after
Lemma 1.3 we arrive at the following fact which is recorded here
for use in the next section.

\medskip

{\bf Lemma 2.1}

{\it The diagram
$$
\xymatrix{
\widehat{G}(K)\ar[r]^-{\subseteq}\ar[d]^{d} &
\bfB_{1}(K)\otimes \Hom_{\dZ}(o,\dZ)\ar[d]^{\log\otimes {\rm id}}
\\
\Hom_{L}(\gfr,K)\ar[r]^-{\subseteq} &
\Hom_{\dQ}(\gfr,K) = K\otimes \Hom_{\dZ}(o,\dZ) }
$$
is cartesian.}

\medskip

We denote by $\Oscr (\widehat{G_0}/K)$ the $K$-Fr\'echet algebra of
global holomorphic functions on the base extension of the variety
$\widehat{G_0}$ to $K$. The main result of Fourier analysis over the
field $\dQ$ is the following.

\medskip

{\bf Theorem 2.2:} (Amice)

{\it The Fourier transform is an isomorphism of $K$-Fr\'echet
algebras}
$$
\matrix{\Fscr : D(G_0 ,K) &
\mathop{\longrightarrow}\limits^{\cong} & \Oscr
(\widehat{G_0}/K)\cr \hfill\lambda & \longmapsto & F_{\lambda}\ .}
$$

Proof: This is a several variable version of [Am2] 1.3 (compare
also [Sc]) based on [Am1].

\medskip

Next we want to compute the ideal $J(o) := \Fscr (I(o))$ in $\Oscr
(\widehat{G_0}/K)$. Let $\xfr_1 := 1 \in \gfr = L$ and $F_t :=
F_{\iota(t\xfr_1)-t\iota (\xfr_1)} \in \Oscr (\widehat{G_0}/K)$
for $t \in L$. A straightforward computation shows that
$$
F_t (\chi_{z \otimes v}) = (tr(tv)-t\cdot tr(v))\cdot\log(z)\ .
$$
By Cor. 1.5 we have the following facts:

1) $\widehat{G}(\dC)$ is the analytic subset of the variety
$\widehat{G_0}/K$ defined by $F_t = 0$ for $t \in L$.

2) $J(o)$ is the ideal of all global holomorphic functions which
vanish on $\widehat{G}(\dC)$.


In 1) one can replace the family of all $F_t$ by finitely many
$F_{t_1},\ldots,F_{t_e}$ if $t_1 ,\ldots,t_e$ runs through a
$\dQ$-basis of $L$.

According to [BGR] 9.5.2 Cor.6 the sheaf of ideals $\Jscr$ in the
structure sheaf $\Oscr_{\widehat{G_0}}$ of the variety
$\widehat{G_0}$ consisting of all germs of functions vanishing on
the analytic subset $\widehat{G}(\dC)$ is coherent. Moreover,
[BGR] 9.5.3 Prop.4 says that the analytic subset
$\widehat{G}(\dC)$ carries the structure of a reduced closed
$L$-analytic subvariety $\widehat{G} \subseteq \widehat{G_0}$ such
that for the structure sheaves we have $\Oscr_{\widehat{G}} =
\Oscr_{\widehat{G_0}}/\Jscr$. Since $\widehat{G_0}$ is a Stein
space the global section functor is exact on coherent sheaves. All
this remains true of course after base extension to $K$. Hence, if
$\Oscr (\widehat{G}/K)$ denotes the ring of global holomorphic
functions on the base extension of the variety $\widehat{G}$ to
$K$ then, by 2), we have
$$
\Oscr (\widehat{G}/K) = \Oscr (\widehat{G_0}/K)/J(o)\ .\leqno(+)
$$
The ideal $J(o)$ being closed $\Oscr (\widehat{G}/K)$ in
particular is in a natural way a $K$-Fr\'echet algebra as well. It
is clear from the open mapping theorem that this quotient topology
on $\Oscr (\widehat{G}/K)$ coincides with the topology as a
projective limit of affinoid algebras.

\eject

{\bf Theorem 2.3:}

{\it The Fourier transform is an isomorphism of $K$-Fr\'echet
algebras}
$$
\matrix{\Fscr : D(G,K) &
\mathop{\longrightarrow}\limits^{\cong} & \Oscr
(\widehat{G}/K)\cr \hfill\lambda & \longmapsto & F_{\lambda}\ .}
$$

Proof: This follows from Thm. 2.2, $(+)$, and the surjection
$(\ast)'$ in section 1.

\medskip



We remark that by construction we (noncanonically) have a
cartesian diagram of rigid $L$-analytic varieties of the form
$$
\xymatrix{
\widehat{G}\ar[r] \ar[d]^{d} &
(\bfB_{1})^{\times d}\ar[d]^{\log}
\\
\Adss^1 \ar[r] &
\Adss^d }
$$
with $\Adss^m$ denoting affine $m$-space where the horizontal arrows
are closed immersions and the vertical arrows are etale. Hence the
variety $\widehat{G}$ is smooth and quasi-Stein (in the sense of
[Kie]).

\bigskip

{\bf 3. Lubin-Tate formal groups and twisted unit disks}

\smallskip

Keeping the notations of the previous section we will give in this
section a different description of the rigid variety $\widehat{G}$. We
will show that the character variety $\widehat{G}$ becomes isomorphic
to the open unit disk after base change to $\dC$.   As a corollary,
the ring of functions ${\cal O}(\widehat{G}/\dC)$ is the same for any
group $G = o$. This result originates in the observation that the
character group $\widehat{G}(\dC)$ can be parametrized with the help
of Lubin-Tate theory.

Fix a prime element $\pi$ of $o$ and let $\Gscr =
\Gscr_{\pi}$ denote the corresponding Lubin-Tate formal group over
$o$. It is commutative and has dimension one and height $[L:\dQ]$.
Most importantly, $\Gscr$ is a formal $o$-module which means that
the ring $o$ acts on $\Gscr$ in such a way that the induced action
of $o$ on the tangent space $t_{\Gscr}$ is the one coming from the
natural $o$-module structure on the latter ([LT]). We always
identify $\Gscr$ with the rigid $L$-analytic open unit disk $\bfB$
around zero in $L$. In this way $\bfB$ becomes an $o$-module
object, and we will denote the action $o \times \bfB
\longrightarrow \bfB$ by $(g,z) \longmapsto [g](z)$. This
identification, of course, also trivializes the tangent space
$t_{\Gscr}$.

Let $\Gscr'$ denote the $p$-divisible group dual to $\Gscr$ and let
$T' = T(\Gscr')$ be the Tate module of $\Gscr'$. Lubin-Tate theory
tells us that $T'$ is a free $o$-module of rank one and that the
Galois action on $T'$ is given by a continuous character $\tau : {\rm
Gal}(\dC/L)
\longrightarrow o^{\times}$. From [Tat] p.177 we know that, by
Cartier duality, $T'$ is naturally identified with the group of
homomorphisms of formal groups over $o_{\dC}$ between $\Gscr$ and the
formal multiplicative group. This gives rise to natural Galois
equivariant and $o$-invariant pairings
$$
<\ ,\,>\ :\;T'\ \mathop{\otimes}\limits_o\ \bfB(\dC)
\longrightarrow
\bfB_1(\dC)
$$
and on tangent spaces
$$
(\ ,\,)\ :\;T'\ \mathop{\otimes}\limits_o\ \dC \longrightarrow
\dC\ .
$$
To describe them explicitly we will denote by $F_{t'}(Z) =
\Omega_{t'}Z + \ldots \in Zo_{\dC}[[Z]]$, for any $t'\in T'$,
the power series giving the corresponding homomorphism of formal
groups. Then
$$
<t',z>\ = 1 + F_{t'}(z)\ \ \ \hbox{and}\ \ \ (t',x) = \Omega_{t'}x\ .
$$

\medskip

{\bf Proposition 3.1:}

{\it The map
$$
\matrix{
\bfB(\dC) \otimes_o T' & \longrightarrow & \widehat{G}(\dC)
\hfill\cr\cr \hfill z \otimes t' & \longmapsto & \kappa_{z\otimes t'}(g)
:= <t',[g](z)> }\leqno{(\diamond)}
$$
is a well defined isomorphism of groups.}

Proof: We will study the following diagram:
$$
\xymatrix{
\bfB(\dC)\otimes_o T'\ar[rr]^(.7){\log_{\Gscr}\otimes {\rm id}}
 \ar@{-->}[rd]\ar[dd]^(.7){\alpha} & &
 t_{\Gscr}(\dC)\otimes_o T' \ar@{-->}[rd]\ar[dd]^(.7){d\alpha}
\\
 & \widehat{G}(\dC) \ar[rr]^(.3){d}\ar[dd]_(.3){\subseteq} &
  & \Hom_{L}(\gfr,\dC)\ar[dd]^(.3){\subseteq}
\\
\Hom_{\dZ}(o,\bfB_{1}(\dC))\ar@{-->}[rd]\ar[rr]^(.7){\Hom(.,\log)\ } & &
\Hom_{\dZ}(o,\dC)\ar@{-->}[rd]
\\
 & \widehat{G}_{0}(\dC)\ar[rr]^(.3){d} &  & \Hom_{\dQ}(\gfr,\dC)
\\
}
$$
Here:{\parindent=10pt
\itemitem{a.} The rear face of the cube in this diagram is the
tensorization by $T'$ of a portion of the map of exact sequences
labelled $(\ast)$ in [Tat] \S4. We use for this the identification
$$
\Hom_{\dZ}(T',.) \otimes_o T' = \Hom_{\dZ}(o,.)\ .
$$
By $\log_{\Gscr}$ we denote the logarithm map of the formal group
$\Gscr$. The map $\alpha$, resp. $d\alpha$, associates to an
element $z\otimes t'$, resp. $\xfr\otimes t'$, the map $g\mapsto
<gt',z>$, resp. $g\mapsto (gt',\xfr)$, for $g\in o$.

\itemitem{b.} The front face of the cube is the diagram
after Lemma 1.3.

\itemitem{c.} The dashed arrows on the bottom face of the cube come
from the discussion before Lemma 2.1.

\itemitem{d.} The dashed arrow in the upper left of the diagram is
the one we want to establish.

\itemitem{e.} The formal $o$-module property of $\Gscr$ says that
the induced $o$-action on $t_{\Gscr}(\dC)$ is the same as the
action by linearity and the inclusion $o \subseteq \dC$. It means
that we have $(gt',\xfr) = (t',g\xfr) = g\cdot(t',\xfr)$ for $g
\in o$ and hence that any map in the image of $d\alpha$ is
$o$-linear. This defines the dashed arrow in the upper right of
the diagram.

}

The back face of the cube is commutative by [Tat] \S4. The front
and bottom faces are commutative by Lemma 1.3 and Lemma 2.1.
Furthermore, the dashed arrows on the bottom of the diagram are
isomorphisms, the left one by the discussion before Lemma 2.1 and
the right one for trivial reasons.

Consider now the right side of the cube. It is commutative by
construction. Since, by [Tat] Prop. 11, $d\alpha$ is injective and
since the lower dashed arrow is bijective the upper dashed arrow
must at least be injective. But by a comparison of dimensions we
see that the dashed arrow in the upper right of the cube is an
isomorphism as well.

In this situation we now may use the fact that, by Lemma 1.3, the
front of the cube is cartesian to obtain that the upper left
dashed arrow is well defined (making the whole cube commutative)
and is given by $(\diamond)$. But according to [Tat] Prop. 11 the
back of the cube also is cartesian. Therefore the map $(\diamond)$
in fact is an isomorphism.\endbeweis

\medskip

Fixing a generator $t_{\rm o}'$ of the $o$-module $T'$ the isomorphism
$(\diamond)$ becomes
$$
\matrix{
\bfB(\dC) & \mathop{\longrightarrow}\limits^{\cong} & \widehat{G}(\dC)
\hfill\cr \hfill z & \longmapsto & \kappa_z := \kappa_{z\otimes
t'_{\rm o}}\ . }\leqno{(\diamond\diamond)}
$$
The main purpose of this section is to see that this latter
isomorphism derives from an isomorphism $\bfB/\dC\,
\mathop{\longrightarrow}\limits^{\cong}\,\widehat{G}/\dC$ between rigid
$\dC$-analytic varieties. In fact, we will exhibit compatible
admissible coverings by affinoid open subsets on both sides.

Let us begin with the left hand side. For any $r \in p^{\Qdss}$ we
have the affinoid disk
$$
\bfB(r) := \{z : |z| \leq r\}
$$
over $L$. Clearly the disks $\bfB(r)$ for $r < 1$ form an
admissible covering of $\bfB$. It actually will be convenient to
normalize the absolute value $|\ |$ and we do this by the
requirement that $|p| = p^{-1}$. The numerical invariants of the
finite extension $L/\dQ$ which will play a role are the
ramification index $e$ and the cardinality $q$ of the residue
class field of $L$. Recall that $o$ acts on $\bfB$ since we
identify $\bfB$ with the formal group $\Gscr_{\pi}$. We need some
information how this covering behaves with respect to the action
of $\pi$.

\medskip

{\bf Lemma 3.2:}

{\it For any $r\in p^{\Qdss}$ such that $p^{-q/e(q-1)} \leq r < 1$ we
have
$$
[\pi]^{-1}(\bfB(r)) = \bfB(r^{1/q})\ \ \ and\ \ \ [p]^{-1}(\bfB(r)) =
\bfB(r^{1/q^e})\ ;
$$
further, in this situation the map
$[\pi^{n}]:\bfB(r^{1/q^{n}})\to\bfB(r)$, for any $n\in \Ndss$, is a
finite etale affinoid map. }

Proof: The second identity is a consequence of the first since $\pi^e$
and $p$ differ by a unit in $o$, and for any unit $g\in o^{\times}$
one has $|[g](z)| = |z|$. Moreover, up to isomorphism, we may assume
([Lan] \S8.1) that the action of the prime element $\pi$ on $\bfB$ is
given by
$$
[\pi](z) = \pi z + z^q\ .
$$
In this case the first identity follows by a straightforward
calculation of absolute values. The finiteness and etaleness of the
map $[\pi^{n}]$ also follows from this explicit formula together with
the fact that a composition of finite etale affinoid maps is finite
etale.

\medskip

Now we consider the right side of $(\diamond\diamond)$. The disk
$\bfB_1$ has the admissible covering by the $L$-affinoid disks
$\bfB_1(r) :=
\{z : |z-1| \leq r\}$ for $r\in p^{\Qdss}$ such that $r < 1$. They are
$\dZ$-submodules so that the $L$-affinoids $\bfB_1(r) \otimes_{\dZ}
\Hom_{\dZ}(o,\dZ)$ form an admissible covering of $\bfB_1
\otimes_{\dZ} \Hom_{\dZ}(o,\dZ)
\cong
\widehat{G_0}$. We therefore have the
admissible covering of $\widehat{G}$ by the $L$-affinoids
$\widehat{G}(r) := \widehat{G} \cap (\bfB_1(r) \otimes_{\dZ}
\Hom_{\dZ}(o,\dZ))$. We emphasize that on both sides the covering
is defined over $L$.

\medskip

{\bf Lemma 3.3:}

{\it For any $r\in p^{\Qdss}$ such that $p^{-p/(p-1)} \leq r < 1$ we
have
$$
\{\chi\in\widehat{G}: \chi^p \in \widehat{G}(r)\} =
\widehat{G}(r^{1/p})\ ;
$$
further, in this situation the map
$[p^{n}]:\widehat{G}(r^{1/p^{n}})\to\widehat{G}(r)$, for any $n\in
\Ndss$, is a finite etale affinoid map. }

Proof: This follows from a corresponding identity between the
affinoids $\bfB_1(r)$. It is, in fact, a special case of the previous
lemma.

\medskip

In order to see in which way the isomorphism $(\diamond\diamond)$
respects these coverings, we need more detailed information on the
power series $F_{t'}$ representing $t'\in T'$.  We summarize the facts
that we require in the following lemma.

\medskip

{\bf Lemma 3.4:}

{\it Suppose $t'\in T'$ is non-zero; the power series
$F_{t'}(Z)=\Omega_{t'}Z+\ldots\in o_{\dC}[[Z]]$ has the following
properties:

a. $\Omega_{gt'}=\Omega_{t'}g$ for $g\in o$;

b.  if $t'$ generates $T'$ as an $o$-module, then
$$
|\Omega_{t'}|=p^{-s}\hbox{\ \ \ with\ \ \ }
s=\frac{1}{p-1}-\frac{1}{e(q-1)}\ ;
$$
c. for any $r<p^{-1/e(q-1)}$, the power series $F_{t'}(Z)$ gives an
analytic isomorphism between $\bfB(r)$ and $\bfB(r|\Omega_{t'}|)$. }

Proof: Part (a) is a restatement of the $o$-linearity of the pairing
$(\ ,\ )$ introduced at the start of this section.  Part (b) follows
from work of Fontaine ([Fon]) on $p$-adic Hodge theory.  We give a
proof in the appendix.  For part (c), recall that $F_{t'}(Z)$ is a
formal group homomorphism.  Therefore if $F_{t'}(z)=0$, $F_{t'}$
vanishes on the entire subgroup of $\bfB(\dC)$ generated by $z$.  The
point $z$ belongs to some $\bfB(r)(\dC)$, and therefore so does the
entire subgroup generated by $z$.  If this subgroup were infinite,
$F_{t'}$ would have infinitely many zeroes in the affinoid
$\bfB(r)(\dC)$ and would therefore be zero.  Consequently $z$ must be
a torsion point of the group $\Gscr$. But other than zero, there are
no torsion points inside the disk $\bfB(r)(\dC)$ if $r<p^{-1/e(q-1)}$
([Lan] \S8.6 Lemma 4 and 5). It follows that the power series
 $F_{t'}(Z)/\Omega_{t'}Z=1+c_1Z+c_2Z^2+\ldots$  has no zeroes inside $\bfB(r)(\dC)$.
Suppose that some coefficient $c_{n}$ in this expansion satisfies
$|c_n|>p^{n/e(q-1)}$.  Then by considering the Newton polygon of the
power series $F_{t'}(Z)/\Omega_{t'}Z$ one sees that the power series
in question must have a zero of absolute value less than
$p^{-1/e(q-1)}$, which we have seen is impossible. Therefore
$|c_{n}|\le p^{n/e(q-1)}$, from which part (c) follows immediately.

\medskip

To simplify the notation, we write $\Omega=\Omega_{t_{\rm o}'}$ for
the ``period'' of the Lubin-Tate group associated to our fixed
generator of $T'$.

By trivializing the tangent space as well as identifying
$\Hom_L(\gfr,\dC)$ with $\dC$ by evaluation at 1 we may simplify the
upper face of the cubical diagram in the proof of Prop. 3.1 to the
following:
$$
\xymatrix{
\bfB(\dC) \ar[rr]^{\log_{\Gscr}}\ar[d]_{z\mapsto \kappa_{z}} & & \dC
\ar[d]^{x\mapsto \Omega x}
\\
\widehat{G}(\dC)\ar[rr]_(.6){\kappa_z\mapsto\log\kappa_z(1)} & & \dC
\\
}\leqno{(\ast)}
$$
Let us examine the map $\kappa(z) := \kappa_{z}$ in coordinates.
Choose for the moment a $\dZ$-basis $e_{1},\ldots,e_{d}$ for $o$
and let $e_{1}^{*},\ldots,e_{d}^{*}$ be the dual basis.  In
coordinates, the map $\kappa:\bfB\to\widehat{G}_{0}$ is given by
$$
\kappa_z=\sum_{i=1}^{d} (1+F_{e_{i}t_{o}'}(z))\otimes e_{i}^{*}\ .
\leqno{(\ast\ast)}
$$
Note first that this map is explicitly rigid $\dC$-analytic and we
know by Prop. 3.1 that this map factorizes through the subvariety
$\widehat{G}/\dC \subset \widehat{G}_{0}/\dC$. We now also see that,
if $r=p^{-q/e(q-1)}<p^{-1/e(q-1)}$, the three parts of Lemma 3.4
together imply that this map carries $\bfB(r)/\dC$ into
$\widehat{G}(r|\Omega|)/\dC$.

\medskip

{\bf Lemma 3.5:}

{\it  Let $r=p^{-q/e(q-1)}$; the map
$$
\matrix{
\bfB(r)/\dC & \mathop{\longrightarrow}\limits^{\cong} &
\widehat{G}(r|\Omega|)/\dC\cr
\hfill z & \longmapsto & \kappa_{z}\hfill\cr }
$$
is a rigid isomorphism. }

Proof: In the discussion preceeding the statement of the lemma we saw
that this  is a well-defined rigid map.  Consider now the other maps
in the diagram $(\ast)$.

-- For $r = p^{-q/e(q-1)}<p^{-1/e(q-1)}$, the logarithm $\log_{\Gscr}$
of the formal group $\Gscr$ restricts to a rigid isomorphism
$$
\log_{\Gscr} : \bfB(r)\ \mathop{\longrightarrow}\limits^{\cong}
\ \bfB(r)\ .
$$
([Lan] \S8.6 Lemma 4).

-- Because $|\Omega|r=p^{-1/(p-1)-1/e}<p^{-1/(p-1)},$ the usual
logarithm restricts to a rigid isomorphism
$$
\log : \bfB_1(r|\Omega|)\ \mathop{\longrightarrow}\limits^{\cong}
\ \bfB(r|\Omega|)\ .
$$
All of this information, together with the diagram $(\ast)$, tells us
that the following diagram of rigid morphisms commutes:
$$
\xymatrix{
\bfB(r)/\dC \ar[rr]^{\cong}\ar[d]_{z\mapsto \kappa_{z}} & &
\bfB(r)/\dC
\ar[d]^{\cong}
\\
\widehat{G}(r|\Omega|)/\dC \ar[rr]_{\kappa_{z}\mapsto\log\kappa_{z}(1)} &
& \bfB(r|\Omega|)/\dC
\\
}
$$

We claim that the lower arrow in this diagram is injective on
$\dC$-points. Assume that $\log\kappa_z(1) = 0$; we then have
$\kappa_z(1) = 1$ which, by the local $L$-analyticity of $\kappa_z$,
means that $\kappa_z$ is locally constant and hence of finite order.
But for our value of $r$ we know that $\bfB_1(r)(\dC)$ is torsionfree
so it follows that $\kappa_z$ must be the trivial character.

Because the upper horizontal and the right vertical map are rigid
isomorphisms and the other two maps at least are injective on
$\dC$-points, all the maps in this diagram  must be isomorphisms on
$\dC$-points.  Because $\widehat{G}$ is reduced, it follows that the
other arrows are rigid isomorphisms as well.


\medskip

This lemma provides a starting point for the proof of the main theorem
of this section.

\medskip

{\bf Theorem 3.6:}

{\it The map
$$
\kappa: \bfB/\dC\ \mathop{\longrightarrow}\limits^{\cong}\ \widehat{G}/\dC
$$
is an isomorphism of rigid varieties over $\dC$; more precisely, if
$r=p^{-q/e(q-1)}$ and $n\in\Ndss_{0}$, then $\kappa$ is a rigid
isomorphism between the affinoids }
$$
\kappa: \bfB(r^{1/q^{en}})/\dC\
\mathop{\longrightarrow}\limits^{\cong}\
\widehat{G}((r|\Omega|)^{1/p^{n}})/\dC\ .
$$
Proof: We remark first that the second statement is in fact stronger
than the first, because  as $n$ runs through $\Ndss_{0}$ the given
affinoids form admissible coverings of $\bfB/\dC$ and
$\widehat{G}/\dC$ respectively.

Lemma 3.5 is the case $n=0$. To obtain the result for all $n$, fix
$n>0$ and consider the diagram:
$$
\xymatrix{
\bfB(r^{1/q^{en}})/\dC \ar[r]^(.44){z\mapsto \kappa_{z}}\ar[d]^{[p^n]} &
\widehat{G}((r|\Omega|)^{1/p^{n}})/\dC \ar[d]^{\chi\mapsto\chi^{p^{n}}}\cr
\bfB(r)/\dC \ar[r]^{z\mapsto\kappa_{z}} & \widehat{G}(r|\Omega|)/\dC \cr
}
$$
By Lemma 3.2, the left-hand vertical arrow is a well-defined finite
etale affinoid map of degree $q^{ne}=p^d$. By Lemma 3.4, part (b), $1
> r|\Omega|=p^{-1/(p-1)-1/e}\ge p^{-p/(p-1)}$ so that Lemma 3.3
applies to the right-hand vertical arrow and it enjoys the same
properties. Lemma 3.5 shows that the lower arrow is a rigid analytic
isomorphism. The upper horizontal arrow then is a well-defined
bijective map on points because of Proposition 3.1 and the first
assertions in Lemma 3.2 and Lemma 3.3. It is a rigid morphism because
it is given in coordinates by the same formula as in the $n=0$ case
(see $(\ast\ast)$).

To complete the argument, let $A$ and $B$ be the affinoid algebras of
$\bfB(r^{1/q^{en}})/\dC$ and $\widehat{G}((r|\Omega|)^{1/p^{n}})/\dC$
respectively. Let $D$ be the affinoid algebra of
$\widehat{G}(r|\Omega|)/\dC$. The rings $A$ and $B$ are finite etale
$D$-algebras of the same rank. The map $f:B\to A$ induced by the upper
arrow in the diagram is a map of $D$-algebras. Because $f$ is
bijective on maximal ideals and $B$ is reduced (because $\widehat{G}$
is reduced), this map is injective. To see that $f$ also is surjective
it suffices, by [B-CA] II\S3.3 Prop. 11, to check that the induced map
$B/\mfr B
\longrightarrow A/\mfr A$ is surjective for any maximal ideal $\mfr \subseteq
D$. But the latter is a map of finite etale algebras over $D/\mfr =
\dC$ of the same dimension which is bijective on points. Hence it
clearly must be bijective.

\medskip

{\bf Corollary 3.7:}

{\it The ring of functions ${\cal O}(\widehat{G}/\dC)$ is isomorphic
to the ring ${\cal O}({\bf B}/\dC)$ of $\dC$-analytic functions on the
open unit disk in $\dC$; in particular, the distribution algebra
$D(G,K)$ is an integral domain.}

\medskip

Let us remark that a careful examination of the proofs shows that
these results in fact hold true over any complete intermediate field
between $L$ and $\dC$ which contains the period $\Omega$.

The ring ${\cal O}({\bf B}/\dC)$ is the ring of power series $F(z) =
\sum_{n\geq 0} a_nz^n$ over $\dC$ which converge on $\{z : |z| < 1\}$.
Let $G_L := {\rm Gal}(\overline{L}/L)$ be the absolute Galois group of
the field $L$ and let $G_L$ act on ${\cal O}({\bf B}/\dC)$ by
$$
F^{\sigma}(z) := \sum_{n\geq 0} \sigma(a_n)z^n\ \ \ \ \hbox{for}\
\sigma\in G_L\ .
$$
By Tate's theorem ([Tat]) we have
$$
\dC^{G_L} = L\ .
$$
Hence the ring ${\cal O}({\bf B})$ coincides with the ring of Galois
fixed elements
$$
{\cal O}({\bf B}) = {\cal O}({\bf B}/\dC)^{G_L}
$$
with respect to this action. This principle which here can be seen
directly on power series in fact holds true for any quasi-separated
rigid $L$-analytic variety $\Xscr$; i.e., one has
$$
\Oscr (\Xscr) = \Oscr (\Xscr /\dC)^{G_L}\ .
$$
By the way the base extension $\Xscr /\dC$ is constructed by pasting
the base extension of affinoids ([BGR] 9.3.6) this identity
immediately is reduced to the case of an affinoid variety. But for any
$L$-affinoid algebra $A$ we may consider an orthonormal base of $A$
and apply Tate's theorem to the coefficients to obtain that
$$
A = (A\ \widehat{\otimes}_L\ \dC)^{G_L}\ .
$$
Since, according to our above theorem, ${\bf B}/\dC$ also is the base
extension of $\widehat{G}$ the ring $\Oscr(\widehat{G})$ must be
isomorphic to the subring of Galois fixed elements in the power series
ring ${\cal O}({\bf B}/\dC)$ with respect to a certain twisted Galois
action. To work this out we first note that the natural Galois action
on $\widehat{G}(\dC)$ is given by composition $\kappa \mapsto
\sigma\circ\kappa$ for $\sigma\in G_L$ and $\kappa\in
\widehat{G}(\dC)$ viewed as a character $\kappa : G \longrightarrow
\dC^{\times}$. Suppose that $\kappa = \kappa_z$ is the image of $z\in
{\bf B}(\dC)$ under the map $(\diamond\diamond)$. The twisted Galois
action $z \mapsto\sigma\ast z$ on ${\bf B}(\dC)$ which we want to
consider then is defined by
$$
\kappa_{\sigma\ast z} = \sigma\circ\kappa_z
$$
and we have
$$
\Oscr(\widehat{G}) \cong \{F\in \Oscr({\bf B}/\dC) : F =
F^{\sigma}(\sigma(\sigma^{-1}\ast z))\ \hbox{for any}\ \sigma\in
G_L\}\ .
$$
Recalling that $\tau : G_L \longrightarrow o^{\times}$ denotes the
character which describes the Galois action on $T'$ we compute
$$
\matrix{
\sigma^{-1}\circ\kappa_{\sigma}(g) & = & <\sigma^{-1}(t_{\rm o}'),
\sigma^{-1}([g](z))>\ =\ < \tau(\sigma^{-1})\cdot t_{\rm
o}',[g](\sigma^{-1}(z))>\hfill\cr\cr & = & <t_{\rm
o}',[g]([\tau(\sigma^{-1})](\sigma^{-1}(z)))>\ =\
\kappa_{[\tau(\sigma^{-1})](\sigma^{-1}(z))}(g)\hfill }
$$
for any $g\in G = o$. Hence
$$
\sigma^{-1}\ast z = [\tau(\sigma^{-1})](\sigma^{-1}(z))\ \ \ \
\hbox{and}\ \ \ \ \sigma(\sigma^{-1}\ast z) = [\tau(\sigma^{-1})](z)\
.
$$

\medskip

{\bf Corollary 3.8:}

$\Oscr(\widehat{G}) \cong \{F\in \Oscr({\bf B}/\dC) : F =
F^{\sigma}\circ [\tau(\sigma^{-1})]\ \hbox{for any}\ \sigma\in G_L\}\
.$

\medskip

The following two negative facts show that our above results cannot be
improved much.

\medskip

{\bf Lemma 3.9:}

{\it Suppose that $K$ is discretely valued; If $L \neq \dQ$ then
$\widehat{G}/K$ and ${\bf B}/K$ are not isomorphic as rigid
$K$-analytic varieties.}

Proof: We consider the difference $\delta_1 - \delta_0 \in D(G,K)$
of the Dirac distributions in the elements $1$ and $0$ in $G = o$,
respectively. Since the image of any character in
$\widehat{G}(\dC)$ lies in the $1$-units of $o_{\dC}$ we see that
the Fourier transform of $\delta_1 - \delta_0$ as a function on
$\widehat{G}(\dC)$ is bounded. On the other hand every torsion
point in $\widehat{G}(\dC)$ corresponding to a locally constant
character on the quotient $o/\dZ$ is a zero of this function.
Since $o \neq \dZ$ we therefore have a nonzero function in
$\Oscr(\widehat{G}/K)$ which is bounded and has infinitely many
zeroes. If $\widehat{G}/K$ and ${\bf B}/K$ were isomorphic this
would imply the existence of a nonzero power series in $\Oscr({\bf
B}/K)$ which as a function on the open unit disk in $\dC$ is
bounded and has infinitely many zeroes. By the maximum principle
the former means that the coefficients of this power series are
bounded in $K$. But according to the Weierstrass preparation
theorem ([B-CA] VII\S3.8 Prop. 6) a nonzero bounded power series
over a discretely valued field can have at most finitely many
zeroes. So we have arrived at a contradiction.

\medskip

According to Lazard ([Laz]) the ring $\Oscr({\bf B}/K)$, for $K$
spherically complete, is a so called Bezout domain which is the
non-noetherian version of a principal ideal domain and which by
definition means that any finitely generated ideal is principal. We
show that this also fails in our setting as soon as $L \neq \dQ$.

\medskip

{\bf Lemma 3.10:}

{\it Suppose that $K$ is discretely valued; if $L \neq \dQ$ then the
ideal of functions in $\Oscr(\widehat{G}/K)$ vanishing in the trivial
character $\kappa_0 \in \widehat{G}(L)$ is finitely generated but not
principal.}

Proof: The ideal in question is a quotient of the corresponding ideal
for the polydisk $\widehat{G_0}/K$ which visibly is finitely
generated. Reasoning by contradiction let $f$ be a generator of the
ideal in the assertion. As a consequence of Theorems A and B ([Kie]
Satz 2.4) for the quasi-Stein variety $\widehat{G}/K$ we then have the
exact sequence of sheaves
$$
0 \longrightarrow \Oscr_{\widehat{G}/K}
\mathop{\longrightarrow}\limits^{f\cdot} \Oscr_{\widehat{G}/K}
\longrightarrow K \longrightarrow 0
$$
on $\widehat{G}/K$ where the third term is a skyscraper sheaf in the
point $\kappa_0$. The corresponding sequence of sections in any
affinoid subdomain is split-exact and hence remains exact after base
extension to $\dC$. It follows that we have a corresponding exact
sequence of sheaves
$$
0 \longrightarrow \Oscr_{\widehat{G}/\dC}
\mathop{\longrightarrow}\limits^{f\cdot} \Oscr_{\widehat{G}/\dC}
\longrightarrow \dC \longrightarrow 0
$$
on $\widehat{G}/\dC$. Using Theorem B we deduce from it that $f$ also
generates the ideal of functions vanishing in $\kappa_0$ in
$\Oscr(\widehat{G}/\dC)$. Consider $f$ now as a rigid map $f :
\widehat{G}/K \longrightarrow \Adss^1/K$ into the affine line. The
composite $f\circ\kappa : {\bf B}/\dC \longrightarrow \widehat{G}/\dC
\longrightarrow \Adss^1/\dC$ then is given by a power series $F\in
\Oscr({\bf B}/\dC)$ which generates the maximal ideal of functions
vanishing in the point $0$. Hence $F$ is of the form
$$
F(z) = az(1 + b_1z + b_2z^2 + \ldots)\ \ \ \ \hbox{with}\ a\in\dC\
\hbox{and}\ b_i\in o_{\dC}
$$
and gives an isomorphism
$$
{\bf B}/\dC \mathop{\longrightarrow}\limits^{\simeq} {\bf
B}^-(|a|)/\dC
$$
between the open unit disk and the open disk ${\bf B}^-(|a|)$ of
radius $|a|$ over $\dC$. We see that $f$ in fact is a rigid map
$$
\widehat{G}/K \longrightarrow {\bf B}^-(|a|)/K
$$
which becomes an isomorphism after base extension to $\dC$. It
follows from the general descent principle we have noted earlier
that $f$ induces an isomorphism of rings
$$
\Oscr({\bf B}^-(|a|)/K) \mathop{\longrightarrow}\limits^{\cong}
\Oscr(\widehat{G}/K)
$$
which obviously respects bounded functions. This leads to a
contradiction by repeating the argument in the proof of the previous
lemma. By a more refined descent argument one can in fact show that
$f$ already is an isomorphism of rigid $K$-analytic varieties which is
in direct contradiction to Lemma 3.9.

\medskip

We close this section by remarking that, because  $\widehat{G}/K$
is a smooth 1-dimensional quasi-Stein rigid variety, one has, for
any $K$, the following positive results about the integral domain
$\Oscr(\widehat{G}/K)$:

1. For ideals $I$ of $\Oscr(\widehat{G}/K)$, the three properties
``$I$ is closed'', ``$I$ is finitely generated'', and ``$I$ is
invertible'' are equivalent.

2. The closed ideals in $\Oscr(\widehat{G}/K)$ are in bijection
with the divisors of $\widehat{G}/K$.  (A divisor is an infinite
sum of closed points, having only finite support in any affinoid
subdomain). In addition a Baire category theory argument shows
that given a divisor $D$ there is a function $F \in
\Oscr(\widehat{G}/K)$ whose divisor is of the form $D + D'$ where
$D$ and $D'$ have disjoint support.

3. Any finitely generated submodule in a finitely generated free
$\Oscr(\widehat{G}/K)$-module is closed.

In particular, $\Oscr(\widehat{G}/K)$ is a Pr\"ufer domain and
consequently a  coherent ring. We omit the proofs which consist of
rather standard applications of Theorems A and B ([Kie]).

\bigskip

{\bf 4. Generalized Mahler expansions}

\smallskip

\def\BB{{\bf B}}
In this section we apply the Fourier theory to obtain a
generalization of the Mahler expansion for locally $L$-analytic
functions. Crucial to our computations is the observation that the
power series $F_{t_{\rm o}'}(Z)$ introduced before Prop. 3.1,
which gives the formal group homomorphism $\Gscr\to\dG$ associated
to $t_{\rm o}'$, is given as a formal power series by the formula
$$
F_{t_{\rm o}'}(Z)=\exp(\Omega\log_{\Gscr}(Z))-1.
$$

Throughout the following, we let $\partial$ denote the invariant
differential on the formal group $\Gscr$.

\medskip

{\bf Definition 4.1:}  {\it For $m\in\Ndss_0$, let $P_{m}(Y)\in
L[Y]$ be the polynomial defined by the formal power series
expansion}
$$
\sum_{m=0}^{\infty} P_{m}(Y)Z^{m}=\exp(Y\log_{\Gscr}(Z))\ .
$$

\medskip

Observe that in the case $\Gscr=\dG$, we have
$$
\exp(Y\log(1+Z))=\sum_{m=0}^{\infty}\left({{Y}\atop{m}}\right)Z^m
$$
so in that case $P_{m}(Y)=\left({{Y}\atop{m}}\right)$.

\medskip

{\bf Lemma 4.2:}

{\it The polynomials $P_{m}(Y)$ satisfy the following properties:
{\parindent=10pt
\itemitem{1.} $P_{0}(Y)=1$ and $P_{1}(Y)=Y$;

\itemitem{2.} $P_{m}(0)=0$ for all $m\ge 1$;

\itemitem{3.}  the degree of $P_{m}$ is exactly $m$, and the leading
coefficient of $P_{m}$ is $1/m!\,$;

\itemitem{4.} $P_{m}(Y+Y')=\sum_{i+j=m}P_{i}(Y)P_{j}(Y')$;

\itemitem{5.} $P_{m}(a\Omega)\in o_{\dC}$ for all $a\in o_{L}$;

\itemitem{6.} for $f(x)\in\dC[[x]]$, we have the identity
$$(P_{m}(\partial)f(x))|_{x=0}=\frac{1}{m!}\frac{d^{m}f}{d
x^{m}}|_{x=0}\ .
$$
} }

Proof: The first four properties are clear from the definition.
The fifth property follows from the fact that, for $a\in o_{L}$,
the power series
$$
F_{at_{\rm o}'}(Z)=F_{t_{\rm
o}'}([a](Z))=\exp(a\Omega\log_{\Gscr}(Z))=\sum_{m=0}^{\infty}
P_{m}(a\Omega)Z^{m}
$$
has coefficients in $o_{\dC}$. For the last property, let
$\delta=\frac{d}{dx}$ be the invariant differential on the
additive formal group.  Then Taylor's formula says that
$$
\exp(\delta b)h(a)=\sum_m (\frac{\delta^{m}}{m!}h(a))b^{m}=h(a+b)\ .
$$
Using the fact that $\log_{\Gscr}$ and $\exp_{\Gscr}$ are inverse
isomorphisms between $\Gscr$ and the additive group over $L$,
Taylor's formula for $\Gscr$ can be obtained by making the
substitutions
$$\matrix{
a & = & \log_{\Gscr}(x)\hfill\cr b & = & \log_{\Gscr}(y)\hfill\cr
h & = & f\circ\exp_{\Gscr}\ . }
$$
It follows easily that $\delta h (a) = \partial f (x)$, and
Taylor's formula becomes
$$
\exp(\partial \log_{\Gscr}(y))f(x)=f(x +_{\Gscr}y)\ .
$$
Comparing coefficients after expanding both sides in $y$ and
setting $x = 0$ gives the result.

\medskip

{\bf Remark:} The identity (6) is part of the theory of Cartier
duality, as sketched for example in Section 1 of [Kat]. Corollary
1.8 of [Kat] shows that $P_{m}(\partial)$ is the invariant
differential operator called there $D(m)$. Comparing this fact
with Formula 1.1 of [Kat] yields the claim in part (6) for the
functions $f(x)=x^n$, and the general fact then follows by
linearity.

\medskip

Our goal now is to study the functions $P_{m}(Y\Omega)$ as
elements of the locally convex vector space $C^{an}(G,\dC)$, where
as always $G=o_L$ as a locally $L$-analytic group. The Banach
space $\Fscr_{a+\pi^no_L}(K)$, for any complete intermediate field
$L \subseteq K \subseteq \dC$ and any coset $a+\pi^no_L$ in $o_L$,
is equipped with the norm
$$
\|\sum_{i=0}^{\infty}
  c_{i}(x-a)^{i}\|_{a,n} := \max_{i}\{|c_{i}\pi^{ni}|\}\ .
$$
This Banach space is the same as the Tate algebra of $K$-valued
rigid analytic functions on the disk $a+\pi^{n}o_{\dC}$, and the
norm, by the maximum principle, has the alternative definition
$$
\|f\|_{a,n}=\max_{x\in a+\pi^{n}o_{\dC}} |f(x)|
$$
which is sometimes more convenient for computation. We recall that
$$
\Fscr_{n}(o_L,K)=\prod_{a\bmod\pi^{n}o_{L}}\Fscr_{a+\pi^no_L}(K)\
\ \ \ \hbox{and}\ \ \ \
C^{an}(G,K)=\lim\limits_{{\longrightarrow}\atop{n}}\,
\Fscr_{n}(o_L,K)\ .
$$

\medskip

{\bf Lemma 4.3:}

{\it For all $a\in o_L$ and all $m\ge 0$, we have
$$
\|P_{m}(Y\Omega)\|_{a,n}\le\max_{0\le i\le m}
\|P_{i}(Y\Omega)\|_{0,n}\ .
$$ }

Proof: By Property (4) of Lemma 4.2, we have
$$
P_{m}((a+\pi^{n}x)\Omega)=\sum_{i+j=m}P_{i}(a\Omega)P_{j}(\pi^{n}\Omega
x)\ .
$$
Therefore, using Property (5) of Lemma 4.2, we have
$$
\matrix{
\|P_{m}(Y\Omega)\|_{a,n}&=\max_{z\in
a+\pi^{n}o_{\dC}}|P_{m}(z\Omega)|\hfill\cr
                        &\le \max_{0\le i\le m}\max_{x\in o_{\dC}}
                        |P_{i}(\pi^{n}\Omega x)|\hfill\cr
            &= \max_{0\le i\le m} \|P_{i}(Y\Omega)\|_{0,n}\
            .\hfill }
$$

\medskip

{\bf Lemma 4.4:}

{\it The following estimate holds for $P_{m}(Y\Omega)$ and $m
\geq 1$:}
$$
||P_{m}(Y\Omega)||_{0,n} <
p^{-1/(p-1)}p^{\frac{m}{eq^{n-1}{(q-1)}}}.
$$

Proof:  Let $r=p^{-q/e(q-1)}$.  As shown in [Lan] \S8.6 Lemma 4,
and as we have used earlier, the functions $\log_{\Gscr}$ and
$\exp_{\Gscr}$ are inverse isomorphisms from $\BB(r)$ to itself,
and $\log_{\Gscr}$ is rigid analytic on all of $\BB$. Furthermore
by Lemma 3.2,  $[\pi^n]^{-1}\BB(r)=\BB(r^{1/q^{n}})$ for any
$n\in\Ndss$. It follows that
$$
\|\log_{\Gscr}(x)\|_{\BB(r^{1/q^{n}})}=rp^{n/e}
$$
where $\|f\|_{\Xscr}$  denotes the spectral semi-norm on an
affinoid $\Xscr$. The function $H(x,y)=y\Omega\log_{\Gscr}(x)$ is
a rigid function of two variables on the affinoid domain
$\BB(r^{1/q^n})\times
\BB(p^{-n/e})$ satisfying
$$
\matrix{
\|H(x,y)\|_{ \BB(r^{1/q^n})\times
\BB(p^{-n/e})} & \leq p^{-n/e}p^{-1/(p-1)+1/e(q-1)}rp^{n/e}\cr\cr
& = p^{-1/(p-1)-1/e} < p^{-1/(p-1)}\ .\hfill }
$$
We conclude from this that $\exp(H(x,y))$ is rigid analytic on $
\BB(r^{1/q^n})\times \BB(p^{-n/e})$ and that
$$
\|\exp(H(x,y))-1\|_{ \BB(r^{1/q^n})\times \BB(p^{-n/e})}<
p^{-1/(p-1)}.
$$
The power series expansion of $\exp(H(x,y))$ is
$$
\exp(H(x,y))=\sum_{m=0}^{\infty} P_{m}(y\Omega)x^{m},
$$
and so we conclude that, for all $m\geq 1$ and for all $y\in
\BB(p^{-n/e})$, we have
$$
|P_{m}(y\Omega)|r^{m/q^n} < p^{-1/(p-1)}.
$$
Therefore we obtain
$$
\|P_{m}(Y\Omega)\|_{0,n} <
p^{-1/(p-1)}p^{\frac{m}{eq^{n-1}(q-1)}}
$$
as desired.

\medskip

From these lemmas we may deduce the following proposition.

\medskip

{\bf Proposition 4.5:}

{\it Given a sequence $\{c_{m}\}_{m\geq 0}$ of elements of $\dC$,
the series
$$
\sum_{m=0}^{\infty} c_{m}P_{m}(y\Omega)
$$
converges to an element of $\Fscr_n(o_L,\dC)$ provided that
$|c_{m}|p^{m/eq^{n-1}(q-1)}\to 0$ as $m\to\infty$. More generally,
this series  converges to an element of $C^{an}(G,\dC)$ provided
that there exists a real number $r$, with $r>1$, such that
$|c_{m}|r^{m}\to 0$ as $m\to\infty$. }

\medskip

Theorems 2.3 and 3.6 together imply the existence of a pairing
$$
\{\ ,\;\} : \Oscr(\BB/\dC)\times C^{an}(G,\dC)\to\dC
$$
that identifies $\Oscr(\BB/\dC)$ and the continuous dual of
$C^{an}(G,\dC)$, both equipped with their projective limit
topologies.  The following lemma gives some basic computational
formulae for this pairing; we will use some of these in the proof
of our main theorem in this section.

\medskip

{\bf Lemma 4.6:}

{\it The following formulae hold for the pairing $\{\ ,\;\}$,
given $F\in\Oscr(\BB/\dC)$ and $f\in C^{an}(G,\dC)$:

{\parindent=10pt
\itemitem{1.} $\{1,f\} = f(0)$;

\itemitem{2.} $\{F_{at_{\rm o}'}, f\} = f(a)-f(0)$ for $a\in o_L$;

\itemitem{3.} $\{F,\kappa_{z}\} = F(z)$ for $z\in\BB$;

\itemitem{4.} $\{F_{at_{\rm o}'}F,f\} = \{F,f(a+.)-f\}$ for $a\in o_L$;

\itemitem{5.} $\{F,\kappa_{z}f\}=\{F(z +_{\Gscr}.),f\}$ for $z\in\BB$;

\itemitem{6.} $\{F,f(a.)\} = \{F\circ [a],f\}$ for $a\in o_L$;

\itemitem{7.} $\{F,f'\} = \{\Omega\log_{\Gscr}\cdot F,f\}$;

\itemitem{8.} $\{F,xf(x)\} = \{\Omega^{-1}\partial F,f\}$;

\itemitem{9.} $\{F,P_{m}(.\Omega)\} = (1/m!)\frac{d^{m}F}{dZ^{m}}(0)$.
}}

Proof: These properties  follow from the definition of the Fourier
transform and from the density of the subspace generated by the
characters in $C^{an}(G,\dC)$.  For example  to see property (5):
if $\lambda'(f):=\lambda(\kappa_{z}f)$, then
$$
F_{\lambda'}(z')=\lambda(\kappa_{z}\kappa_{z'})=
\lambda(\kappa_{z +_{\Gscr} z'})=F_{\lambda}(z +_{\Gscr} z')\ .
$$
For property (7), using (4) we have
$$
\matrix{
\hfill \{F,f'\}&=\mathop{\lim}\limits_{\epsilon\to 0}\epsilon^{-1}
\{F,f(.+\epsilon)-f\}\hfill \cr\cr
&=\mathop{\lim}\limits_{\epsilon\to 0}\epsilon^{-1}\{F_{\epsilon
t_{\rm o}'}F,f\}\hfill
\cr\cr &=\{\Omega\log_{\Gscr}\cdot F,f\}\hfill }
$$
using continuity and the fact that $F_{\epsilon t_{\rm
o}'}=\exp(\epsilon\Omega\log_{\Gscr}) - 1$ for small $\epsilon$.
An analogous computation based on (5) gives (8). The last property
(9) follows from Lemma 4.2.6 and (8).

\medskip

We may now prove the main result of this section.

\eject

{\bf Theorem 4.7:}

{\it Any function $f\in C^{an}(G,\dC)$ has a unique representation
in the form
$$
f=\sum_{m=0}^{\infty} c_m P_{m}(.\Omega)
$$
as in Prop. 4.5; in this representation,   $c_m=\{Z^{m},f\}$. }

Proof: Part (9) of Lemma 4.6, along with continuity, shows that if
$f$ has a representation in the given form then $c_m=\{Z^{m},f\}$.
The functions $Z^{m}$ generate a dense subspace in
$\Oscr(\widehat{G}/\dC)$, and so a function $f$ with all $c_m=0$
must be zero (for example all Dirac distributions pair to zero
against $f$).  This implies that this type of  representation, if
it exists, is unique. Suppose we show that, for any $f\in
C^{an}(G,\dC)$, there exists an $r>1$ such that
$|\{Z^{m},f\}|\,r^{m}\to 0$ as $m\to\infty$.  Then by Prop. 4.5,
the series
$$
\overline{f}(x):=\sum_{m=0}^{\infty} \{Z^{m},f\}P_{m}(x\Omega)
$$
converges to a locally analytic function and by Lemma 4.6, Part 9,
we have $\{Z^{m},f\}=\{Z^{m},\overline{f}\}$ for all $m$.
Therefore $\overline{f}=f$.

Thus we have reduced our main theorem to the claim that
$|\{Z^{m},f\}|\,r^{m}\to 0$ as $m\to \infty$ for some $r>1$. The
function $f$ being locally analytic it belongs to one of the
Banach spaces $\Fscr_{n}(o_L,\dC)$. Using the topological
isomorphism between the Fr\'{e}chet spaces $\Oscr(\BB/\dC)$ and
$D(G,\dC)=C^{an}(G,\dC)'$, there is a rational number  $s>0$ such
that the map $\Oscr(\BB/\dC)\to\Fscr_{n}(o_L,\dC)'$ factors
through the Tate algebra $\Oscr(\BB(p^{-s})/\dC)$.  If we choose
another rational number $s'$ so that  $0<s'<s$, then in the Tate
algebra $\Oscr(\BB(p^{-s})/\dC)$, the set of rigid functions
$\{(Z/p^{-s'})^{m}\}_{m\geq 0}$ goes to zero and therefore so does
$|\{Z^{m},f\}|\,p^{s'm}$. This proves the existence of the desired
expansion.

\medskip

{\bf Remark:} In [Kat], Katz discusses what he calls
``Gal-continuous'' functions.  These are continuous functions on
$G$ that satisfy (in our notation) $\sigma(f(x))=f(\tau(\sigma)x)$
for all $\sigma\in G_{L}$. If $\{c_{m}\}$ is a sequence of
elements of $L$ such that $|c_m|\to 0$, then
$f(x):=\sum_{m=0}^{\infty} c_{m}P_{m}(x\Omega)$ is continuous by
Part (5) of Lemma 4.2, and by the Galois properties of $\Omega$ it
is even Gal-continuous.

\bigskip

{\bf 5. $p$-adic L-functions}

\smallskip

In this section we will illustrate how the integration theory
developed in this paper applies to yield $p$-adic L-functions for
CM elliptic curves $E$ at supersingular primes. In fact, our
method allows us to apply the Coleman power series approach
described in [dS] directly in the supersingular case. We will
content ourselves with proving a supersingular analogue of a weak
version of Theorem II.4.11 of [dS]; this should demonstrate
sufficiently the nature of our construction, without requiring too
much of a diversion into global arithmetic.

We should emphasize that the L-functions we will construct in the
supersingular case come from locally analytic distributions on
Galois groups, not measures.  A character on the Galois group is
integrable provided that its restriction to a small open subgroup
is a power of $\ophi$, where $\ophi$ gives the representation on
the dual Tate module of $E$. We therefore cannot make any
immediate connection to Iwasawa module structure of, for example,
elliptic units. This is an interesting problem for the future.

Our results are closely related to those of Boxall ([Box]). See
the Remark after Prop. 5.2 for more discussion of the
relationship.

Before discussing $p$-adic L-functions we will develop Fourier
theory for the multiplicative group;   this will be useful because
the $p$-adic L-functions we construct arise as locally analytic
distributions on Galois groups that are naturally isomorphic to
multiplicative, rather than additive groups.  Let $H$ be
$o_L^{\times}$ as $L$-analytic group and let $H_1$ be the subgroup
$1+\pi o_L$. Using the Teichm{\"u}ller character $\omega$, we have
$$
H=H_1\times k^{\times}
$$
where $k$ is the residue field of $o_L$. For $x\in H$, let $<x>$
be the projection of $x$ to $H_1$.

As always, $G$ is the additive group $o_L$. Let us assume that the
absolute ramification index $e$ of the field $L$ satisfies
$e<p-1$. We define $\ell := \pi^{-1}\cdot\log$ so that $\ell: H_1
\mathop{\longrightarrow}\limits^{\cong} G$ is an $L$-analytic isomorphism.

This map induces an isomorphism between the distribution algebras
$D(H_1,K)$ and $D(G,K)$.   The group $\widehat{H}(\dC)$ of locally
$L$-analytic, $\dC$-valued characters of $H$ is isomorphic to a
product of $q-1$ copies of the open unit disk $\bfB$ using the
results of section 3, indexed by the (finite) character group of
$k^{\times}$.  For $z\in\bfB(\dC)$, let $\psi_{z}$ be the
corresponding character of $H_{1}$.  Then for any distribution
$\lambda\in D(H,K)$, and any character $\omega^{i}\psi_{z}$ with
$z\in\bfB(\dC)$, and $0\le i\le q-1$, we have the ``Mellin
transform''
$$
M_{\lambda}(z,\omega^{i})=\lambda(\omega^{i}\psi_{z})\ .
$$
For each fixed value of the second variable, $M_{\lambda}$ is a
rigid function in $\Oscr(\bfB/\dC)$.

Now let us compare the Fourier transforms for $G$ and $H$ in a
different way.  The group $o_L^{\times}$, as an $L$-analytic
manifold, is an open submanifold of $o_L$.  If we have a
distribution $\lambda$ in $D(G,K)$ that vanishes on functions with
support in $\pi o_L$, then $\lambda$ gives a distribution on
$H=o_L^{\times}\subset o_L$.  It follows easily from Lemma 4.6.5
that $\lambda$ is supported on $H$ precisely when its Fourier
transform $F_{\lambda}$ satisfies
$$
\sum_{[\pi](z)=0} F_{\lambda}(.+_{\Gscr}z)=0\ .
$$
We have the following result comparing the Fourier and Mellin
transforms.

\medskip

{\bf Proposition 5.1:}

{\it Let $\lambda$ be a distribution in $D(G,\dC)$ supported on
$H$, let $F_{\lambda}$ be its Fourier transform, and let
$M_{\lambda}$ be its Mellin transform; suppose that $n\in\Ndss$
satisfies $n\equiv i\pmod{q-1}$; then}
$$
M_{\lambda}(\exp_{\Gscr}(n\pi/\Omega),\omega^{i})=
\mathop{\int}\limits_{o^{\times}_L} x^n d\lambda(x)
= \Omega^{-n}(\partial^{n}F_{\lambda}(z)|_{z=0})\ .
$$

Note that the hypothesis $e<p-1$ guarantees that
$M_{\lambda}(\exp_{\Gscr}(x\pi/\Omega),\omega^i)$ is a (globally)
analytic function of $x\in o_L$.  Thus the left hand side of these
equations gives a (globally) $L$-analytic interpolation of the
values on the right side.

Proof: Let $z(n)=\exp_{\Gscr}(n\pi/\Omega)$.
 By definition,
$$
M_{\lambda}(z(n),\omega^{i})=\lambda(\omega^{i}\psi_{z(n)})\ .
$$
Now
$$
\matrix{
\psi_{z(n)}(<x>)&=&\kappa_{z(n)}(\ell(<x>))\hfill\cr
&=& t'_{\rm o}([\ell(<x>)](z(n)))\hfill\cr &=&
\exp(\Omega\ell(<x>)\log_{\Gscr}(z(n)))\hfill  }
$$
because $|[\ell(<x>)](z(n))|<p^{-1/e(q-1)}$ (by Lemma 3.4.b and
the hypothesis $e<p-1$). But
$$
\exp(\Omega\ell(<x>)\log_{\Gscr}(z(n)))
=\exp(n\log(<x>))=<x>^{n}
$$
so $(\omega^{i}\psi_{z(n)})(x)=x^{n}$.  But
$$
\lambda(x \mapsto x^{n})=\Omega^{-n}(\partial^{n}F_{\lambda}(z)|_{z=0})
$$
by Lemma 4.6.8/9.

\medskip

Now we will embark on a digression into the theory of CM elliptic
curves, following the notation and the logic of Chap. II in [dS].
Let $\bfK$ be an imaginary quadratic field, and let $\ffr$ be an
integral ideal of $\bfK$ such that the roots of unity in $\bfK$
are distinct mod $\ffr$. Let $p$ be a rational prime that is
relatively prime to $6\ffr$ and inert in $\bfK$. Let $\bfF$ be the
ray class field $\bfK(\ffr)$ and let $\bfF_{n}:=\bfK(p^{n}\ffr)$
and $\bfF_{\infty} := \bigcup_{n\in\Ndss} \bfF_{n}$. Assume for
technical reasons that will become clear in a moment that $p$ as a
prime of $\bfK$ {\it splits completely in} $\bfF$. Let $\wp$ be a
prime above $p$ in $\bfF$. The prime $\wp$ ramifies totally in
$\bfF_{\infty}$; let $F_{\infty}$ be the completion of
$\bfF_{\infty}$ at the unique prime above $\wp$. Let $o$ be the
ring of integers in the local field $\bfK_{p}$.

Fix an elliptic curve $E$ over $\bfF$ with CM by the ring of
integers in $\bfK$ and with associated Hecke character of the form
$\psi_{E/\bfF}=\varphi\circ N_{\bfF/\bfK}$, where $\varphi$ is a
Hecke character of $\bfK$ of type $(1,0)$ and conductor dividing
$\ffr$; we moreover assume that there is a complex period
$\Omega_{\infty}$ so that the period lattice $\Lscr$ of $E$ is
$\Omega_{\infty}\ffr$. We view $\varphi$ also as a
$\bfK_{p}^{\times}$-valued character of
$\Gamma_{\bfK}:=Gal(\bfF_{\infty}/\bfK)$. If $\afr$ is an integral
ideal of $\bfK$ such that the Artin symbol $\sigma_{\afr}$ belongs
to $\Gamma_{\bfF}=Gal(\bfF_{\infty}/\bfF)$, then $\sigma_{\afr}$
acts on the $p$-adic Tate module of $E$ through multiplication by
$\varphi(\afr)$. We let $\ophi$ be the Hecke character giving the
action of $\Gamma_{\bfF}$ on the dual Tate module of $E$.  The
character $\ophi$ gives us an isomorphism
$$
\ophi:\Gamma_{\bfF}\to o^{\times}.
$$
We use this isomorphism to equip $\Gamma_{\bfF}$ with an
$o$-analytic structure. We let ${\bf N}$ denote the absolute norm.

Our assumption that $p$ splits completely in $\bfF$ means that
that the formal group $\widehat{E}_{\wp}$ of $E$ at $\wp$ is a
Lubin-Tate group over $o$ of height two.  (To handle general $p$,
deShalit works with what he calls ``relative'' Lubin-Tate groups.
Presumably one can generalize our results to this situation as
well.)   Furthermore, the field $\bfF_{\infty}$ contains all of
the $p$-power torsion points of this formal group, as well as (by
the Weil pairing)  all of the $p$-power roots of unity.  Thus our
uniformization result holds over this field. Choose an
$o$-generator $t'_{\rm o}$ of the (global) dual Tate module
$\Hom(T_{p}(E),T_p(\dG))$ defined over $\bfF_{\infty}$. Then  the
pairing $\{\ ,\,\}$ from section 4 looks like:
$$
\Oscr(\widehat{E}_{\wp}/F_{\infty})\times C^{an}(o,F_{\infty})\to F_{\infty}\ .
$$

Now we introduce the machinery of Coleman power series and
elliptic units. Let $\afr$ be an integral ideal relatively prime
to $p\ffr$ and let $\Theta(y;\Lscr,\afr)$ be the elliptic function
from [dS] II.2.3 (10). Let
$Q_{\afr}(Z)=\Theta(\Omega_{\infty}-\log_{\widehat{E}_{\wp}}(Z);\Lscr,\afr)$
in $o[[Z]]$ (see [dS] Prop. II.4.9; note that this proposition is
true for inert primes as well as split ones, as is clear from its
proof). The power series $Q_{\afr}(Z)$ is the Coleman power series
associated to a norm-compatible sequence of elliptic units, as
deShalit explains.

Define
$$
g_{\afr}(Z)=\log Q_{\afr}(Z)-{{1}\over{p^{2}}}\sum_{\mathop{z\in
\widehat{E}_{\wp}}\limits_{[p](z)=0}}\log Q_{\afr}(Z +_{\widehat{E}_{\wp}} z)\ .
$$

\medskip

{\bf Proposition 5.2:}

{\it The power series
$g_{\afr}(Z)\in\Oscr(\widehat{E}_{\wp}/F_{\infty})$ is the Fourier
transform of an $F_{\infty}$-valued, locally analytic distribution
on $o$ supported on $o^{\times}$.  By means of the isomorphism
$\ophi$ from $\Gamma_{\bfF}$ to $o^{\times}$, it defines a locally
analytic distribution on $\Gamma_{\bfF}$ with the interpolation
property
$$
\matrix{
\Omega^{m}\{g_{\afr}(Z),\ophi^{m}\}\hfill\cr\cr
   \ \ \ =-12(1-p^{m-2})\Omega_{\infty}^{-m}(m-1)!({\bf N}(\afr)
   L_{\ffr}(\ophi,m,1)-\varphi(\afr)^{m}L_{\ffr}(\ophi,m,\afr)) }
$$
for any $m\in\Ndss$. Here $L_{\ffr}(\ophi,s,\cfr)$ denotes the
``partial'' Hecke L-function of conductor $\ffr$, equal to $\sum
_{\bfr}\ophi(\bfr){\bf N}(\bfr)^{-s}$ over ideals $\bfr$ prime to
$\ffr$ and such that $(\bfr,\bfF/\bfK)=(\cfr,\bfF/\bfK)$.}

Proof: The first assertion follows easily from the formulae in
Lemma 4.6.  The interpolation property comes from the formula
(Lemma 4.6 again):
$$
\{g_{\afr}(Z),\ophi^{m}\}=\Omega^{-m}\partial^{m}g_{\afr}(Z)|_{Z=0}\ .
$$
The rest of the computation is just a version of [dS] II.4.10. The
invariant differential $\partial$ pulls back to $d/dy$ on the
complex uniformization of $E$, so
$$
\Omega^{m}\{g_{\afr}(Z),\ophi^{m}\}=\left({{d}\over{dy}}\right)^{m}(\log
\Theta(\Omega_{\infty}-y;\Lscr,\afr) - p^{-2}\log
\Theta(\Omega_{\infty}-y;p^{-1}\Lscr,\afr))|_{y=0}
$$
and the claimed formula then follows from the equivalent in our
situation of [dS] II.4.7 (17), along with  II.3.1 (7) and Prop.
II.3.5.

\medskip

{\bf Remark:} From Prop. 5.1, we see that the Mellin transform
$M_{\afr}$ of the distribution $\{g_{\afr}(Z),\cdot\}$ is an
$o$-analytic function on $o$ interpolating the special values
$\Omega^{-m}\partial^{m}g_{\afr}(Z)|_{Z=0}$. This function (on
$\dZ$) was  constructed by Boxall ([Box]). Other than the fact
that our construction is arguably more natural, the principal new
results here are that our function is $\bfK_{p}$-analytic on $o$
rather than $\dQ$-analytic on $\dZ$. The existence of such an
analytic interpolating function implies congruences among the
special values.

\medskip

We may give the slightly larger Galois group $\Gamma_{\bfK}$ a
locally  analytic structure by transporting that of
$\Gamma_{\bfF}$ to its finitely many cosets in $\Gamma_{\bfK}$. To
extend the integration pairing to $\Gamma_{\bfK}$, recall that,
along with $E$ we have finitely many other elliptic curves
$E^{\sigma}$ as $\sigma$ runs through $Gal(\bfF/\bfK)$. We also
have, for each $\sigma\in\Gamma_{\bfK}$, an isogeny
$\iota(\sigma):E \to E^{\sigma}$  as in [dS] Prop. II.1.5. If
$\afr$ is an ideal prime to $p$, then the associated isogeny
$\iota(\sigma_{\afr})$ has degree ${\bf N}(\afr)$, which is prime
to $p$ and therefore induces an isomorphism between the formal
groups $\widehat{E}_{\wp}$ and $\widehat{E^{\sigma_{\afr}}}_{\wp}$
.

A typical locally analytic function $\overline{f}$  on
$\Gamma_{\bfK}$ may be written
$$
\overline{f}(\sigma)=f_{i}(\ophi(\sigma_{i}^{-1}\sigma))\hbox{\rm\quad when
$\sigma\in\sigma_{i}\Gamma_{\bfF}$}
$$
where $\afr_{i}$ is a collection of integral ideals of $\bfK$ so
that the Artin symbols $\sigma_{i}=\sigma_{\afr_{i}}$ form the set
$Gal(\bfF/\bfK)$, and $f_{i}\in C^{an}(o,F_{\infty})$, supported
on $o^{\times}$.

Let $\Oscr(\widehat{E}_{\wp}/F_{\infty})^{0}$ denote the subspace
of functions $F \in \Oscr(\widehat{E}_{\wp}/F_{\infty})$
satisfying\break $\sum_{[p](z)=0} F(. +_{\widehat{E}_{\wp}} z)=0$.
These are the distributions supported on $o^{\times}$. We define
an extended integration pairing
$$
\{\ ,\,\} : \mathop{\oplus}_{\sigma_{\afr}\in Gal(\bfF/\bfK)}
\Oscr(\widehat{E^{\sigma_{\afr}}}_{\wp}/F_{\infty})^{0}
\times C^{an}(\Gamma_{\bfK},F_{\infty}) \longrightarrow F_{\infty}
$$
by setting
$$
\{\overline{h},\overline{f}\} := \sum_{i}\{h_{i}\circ\iota(\sigma_{i}),f_{i}
({\bf N}(\afr_{i})^{-1}.)\}\ .
$$

\medskip

{\bf Lemma 5.3:}

{\it This pairing is well-defined (i.e., it is independent of the
choice of coset representatives), and identifies the left hand
space with the continuous dual of the right hand space. }

Proof:  The duality is clear; the key point is that the pairing is
well defined. Suppose we replace $\sigma_{i}$ with
$\sigma_{i}\tau_{\bfr}$, where $\tau_{\bfr}\in\Gamma_{\bfF}$. Then
$h_i\circ\iota(\sigma_{i}\tau_{\bfr})=h_i\circ\iota(\sigma_{i})\circ[\varphi
(\tau_{\bfr})]$ (see [dS] II.4.5).  The decomposition of
$\overline{f}$ also changes, with $f_{i}$  replaced by
$f'_{i}=f_{i}(\ophi(\tau_{\bfr}).)$.  Then, using Lemma 4.6 as
usual, the pairing satisfies
$$
\{h_{i}\circ\iota(\sigma_{i})\circ[\varphi(\tau_{\bfr})],f_{i}
({\bf N}(\afr_{i})^{-1}{\bf N}(\bfr)^{-1}\ophi(\tau_{\bfr}).)\}
=\{h_{i}\circ\iota(\sigma_{i}),f_{i}({\bf N}(\afr_{i})^{-1}.)\}\ .
$$

\medskip

{\bf Theorem 5.4:} {\it (Compare [dS] Thm. II.4.11)

Let $\overline{h}=\{h_{i}\}$ where $h_{i}:=\sigma_{i}(g_{\afr})$
with $g_{\afr}$ the (formal) elliptic function over $\bfF$ used
for the construction of the partial L-function in Prop. 5.2. Let
$\epsilon$ be any locally analytic character on $\Gamma_{\bfK}$,
whose restriction to $\Gamma_{\bfF}$ is $\ophi^{m}$ for some
$m\in\Ndss$. Then the locally analytic distribution $\overline{h}$
on $\Gamma_{\bfK}$ has the interpolation property }
$$
\Omega^{m}\{\overline{h},\epsilon\}=12(m-1)!\Omega_{\infty}^{-m}
(1-N(p)^m\epsilon(p)^{-1}p^{-2}) ({\bf
N}(\afr)^{m}\epsilon(\afr)^{-1}-{\bf N}
(\afr))L_{\ffr}(\epsilon,m)\ .
$$
Proof: The proof is a long computation very much in the spirit of
[dS] Thm. II.4.11 (though we have cheated in the statement of the
Theorem and avoided the case where $p$ divides the conductor).
Choose coset representatives $\sigma_{i}=\sigma_{\cfr_{i}}$. The
point of the computation is that ([dS] II.2.4 (ii) and II.4.5
(iv))
$$
\sigma_{i}(g_{\afr})\circ\iota(\sigma_{i})=g_{\afr\cfr_{i}}-{\bf N}
(\afr)g_{\cfr_{i}}
$$
and
$$
f_{i}(a) = \epsilon(\sigma_{i})a^{m}\ .
$$
Then the partial terms in the pairing are
$$
{\bf
N}(\cfr_{i})^{-m}\epsilon(\sigma_{i})\Omega^{-m}(\partial^{m}g_{\afr\cfr_{i}}-{\bf
N}(\afr)\partial^{m}g_{\cfr_{i}})\ .
$$
These terms may then be evaluated using Prop. 5.2, and when the
results are combined one obtains the statement of the theorem.

\medskip

{\bf Remark:} One can compute an interpolation result for more
general locally analytic characters $\epsilon$ -- explicitly,
characters which restrict to $\ophi^{m}$ on an open subgroup of
$\Gamma_{\bfF}$ --  by following the same line of argument as in
[dS] II.4.11.

\bigskip

{\bf Appendix. $p$-adic periods of Lubin-Tate groups}

\smallskip

In the analysis in section 3 of the behavior of the isomorphism
$(\diamond\diamond)$ relative to the affinoid coverings  on
$\widehat{G}$ and on $\bfB$ we needed rather exact information about
the ``period'' $\Omega$ of the Lubin-Tate group $\Gscr$.  In this
appendix, we apply results of Fontaine [Fon] to obtain this
information.

All of the significant ideas in this section come from the article
[Fon], and we follow the notation of that article with the
following exceptions. We will use the letter $X$ for the module of
differentials $\Omega_{o_{L}}(o_{\overline{L}})$ (called $\Omega$
by Fontaine). We also do not distinguish between $\Gscr$ as a
formal group or $p$-divisible group thanks to [Tat] Prop. 1.

As before, we let $\Gscr$ be the Lubin-Tate group over $o$ associated
to the uniformizer $\pi$ and let $\Gscr'$ be the dual $p$-divisible
group. We denote by $q$ and $e$ be respectively the number of elements
in $o/\pi o$ and the ramification index of $L/\dQ$. Furthermore, $T$
and $T'$ are the Tate modules of $\Gscr$ and $\Gscr'$ respectively,
and $F_{t'}(Z)\in Zo_{\dC}[[Z]]$ is the power series corresponding to
$t'\in T'$ as in section 3. We let $\omega$ be the invariant
differential on $\Gscr$ such that $\omega=(1+\ldots)dZ$. We define
$\Omega_{t'}$ so that $F_{t'}(Z)=\Omega_{t'}Z+\cdots.$

We write $\Gscr_{n}$ and $\Gscr_{n}'$ for the group schemes of $p^{n}$
torsion points on $\Gscr$ and $\Gscr'$ respectively.  Let $L_{n}$ be
the finite extension field of $L$ generated by the
$\overline{L}$-points of $\Gscr_{n}'$.

The various maps that are denoted by decorated forms of the letter
$\phi$ are those defined in [Fon].

We remark that, in the case that $L/\dQ$ is tamely ramified, so that
$e\le p-1$, the result of part (c) of the following Theorem was
obtained by Boxall ([Box]) by power series computations.

{\bf Theorem:}

{\it a. For $t'\in T'$, we have
$$
\phi^{0}_{\Gscr'}(t')=\frac{dF_{t'}}{1+F_{t'}}=\Omega_{t'}\omega\ ;
$$
b. there is a sequence of elements $\Omega_{t'}(n)\in o_{L_{n}}$ for
integers $n\ge 1$ such that $\Omega_{t'}(n+1)\equiv
\Omega_{t'}(n)\pmod{p^{n}o_{L_{n+1}}}$ and such that the sequence of
$\Omega_{t'}(n)$ converges in $\dC$ to $\Omega_{t'}\,$;

c. let $t'_{o}$ be any generator of the $o$ module $T'\,$; then the
fundamental period $\Omega=\Omega_{t'_{o}}$ satisfies
$$
|\Omega|=p^{-s}
$$
where}
$$
s=\frac{1}{p-1}-\frac{1}{e(q-1)}\ .
$$

Proof: We begin with parts (a) and (b).  First of all, $F_{t'}(Z)$
being a formal group homomorphism from $\Gscr$ to $\dG$, the pullback
of the invariant differential $dZ/(1+Z)$ on $\dG$ must be a multiple
of $\omega$.  Comparing leading coefficients shows that
$dF_{t'}(Z)/(1+F_{t'}(Z))=\Omega_{t'}\omega$.  Fontaine's map
$$
\phi^{0}_{\Gscr'}: o_{\dC}\otimes_{\dZ}T'\to t^{*}_{\Gscr}(o_{\dC})
$$
as defined on p. 406 of [Fon], is the limit of maps
$$
\phi^{0}_{\Gscr'_{n}}:o_{\overline{L}}\otimes\Gscr'_{n}(o_{\overline{L}})\to
t^{*}_{\Gscr_{n}}(o_{\overline{L}})
$$
defined on p. 396.    To compute these maps, recall that the affine
algebra of $\Gscr_{n}$ is $R_{n}=o_{L}[[Z]]/J_{n}$ where $J_{n}$ is
the ideal generated by $[p^{n}](Z)$.  The element $t'$ is represented
explicitly as $(t'(n))_{n}$ where the $t'(n)$ are a compatible
sequence of homomorphisms $\Gscr_{n}\to\mu_{p^{n}}$ over $o_{\dC}$.
Further, the element $t'(n)$ is given explicitly by the class
$1+F_{t'}(Z)+J_{n}$ in $o_{\dC}\otimes R_{n}$. But each homomorphism
$\Gscr_{n}\to\mu_{p^{n}}$ is defined over $o_{L_{n}}$. Therefore
$1+F_{t'}(Z)\equiv g_{n}(Z)\pmod{J_{n}}$ for some $g_{n}\in
o_{L_{n}}\otimes R_{n}$. Now on the one hand Fontaine's map is given
by the formula
$$
\phi^{0}_{\Gscr_{n}'}(t'(n))=\frac{dg_{n}}{g_{n}}\in
t^{*}_{\Gscr_{n}}(o_{L_{n}})\ .
$$
By Prop. 10 of [Fon] we know that
$$
t^{*}_{\Gscr_{n}}(o_{L_{n}})=t^{*}_{\Gscr}(o_{L_{n}})/p^{n}t^{*}_{\Gscr}(o_{L_{n}})\
.
$$
Therefore we define $\Omega_{t'}(n)\in o_{L_{n}}$ so that
$$
dg_{n}/g_{n}\equiv
\Omega_{t'}(n)\omega\pmod{p^{n} t^{*}_{\Gscr}(o_{L_{n}})}\ .
$$
But both $g_{n}$ and $1+F_{t'}$ represent the same map over $o_{\dC}$
from $\Gscr_{n}$ to $\mu_{p^{n}}$, and therefore we must have
$$
\Omega_{t'}(n)\equiv\Omega_{t'}\pmod{p^{n} t^{*}_{\Gscr}(o_{\dC})}\ .
$$
By definition $\phi^{0}_{\Gscr'}(t')$ is the limit of the
$\Omega_{t'}(n)\omega$, which we have just shown is
$\Omega_{t'}\omega\,$.

Now consider the following commutative diagram, obtained from Prop. 8
of [Fon] by applying section 5.9 to pass to the inverse limit over
multiplication by $p\,$:
$$
\xymatrix{ o_{\dC}\otimes_{\dZ} T\ \times\ o_{\dC}
\otimes_{\dZ} T'
 \ar[r]^(.37){\phi}\ar[d]^{\theta}
      &
 t_{\Gscr'}^{*}(o_{\dC})
\oplus t_{\Gscr}(T_{p}(X))
\ \times\
t^{*}_{\Gscr}(o_{\dC})
\oplus t_{\Gscr'}(T_{p}(X))\ar[d]^{\nu}
\\
o_{\dC}\otimes T_{p}(\dG) \ar[r]^{\xi}   & T_{p}(X) \\ }
$$
Here the map $\phi$ is
$$
\phi={\phi_{\Gscr,o_{\dC}}\times \phi_{\Gscr',o_{\dC}}}
$$
as defined in Prop. 11 of [Fon], where it is shown to be injective,
and to induce an isomorphism upon tensoring with $\dC$. The vertical
arrows are the natural pairings, and the lower horizontal arrow is
induced by the map $\xi$ of Thm. 1' of [Fon].

Each of the spaces $\dC\otimes_{\dZ} T$ and $\dC\otimes_{\dZ}T'$
decompose into a direct sum of  one-dimensional eigenspaces
corresponding to distinct embeddings of $L\into\dC$.  The map $\phi$
is $o$-linear and therefore must respect this decomposition.  On the
upper right, the $o$-actions on the spaces $t^{*}_{\Gscr}(o_{\dC})$
and $t_{\Gscr}(T_{p}(X))$ are given by the given embedding $o
\subseteq o_{\dC}$. Therefore the above diagram can be ``reduced'' to
the following:

$$
\xymatrix{
(o_{\dC}\otimes_{o} T)\times \Hom_o(T,o_{\dC}(1))
\ar[r]^(.56){\phi}\ar[d]^{\theta}
      &  t_{\Gscr}(T_{p}(X))\times t^{*}_{\Gscr}(o_{\dC})\ar[d]^{\nu}\\
o_{\dC}\otimes T_{p}(\dG) \ar[r]^{\xi}   &T_{p}(X) \\ }
$$
Choose  a generator $u$ of $T$, a generator $\epsilon$ of
$T_{p}(\dG)$, and let $f\in \Hom_o(T,o_{\dC}(1))$ be the unique
$o$-linear map such that $f(u)=\epsilon$.  Trace the pairing $(u,f)$
both ways through the square, using the explicit formulae for the maps
involved from [Fon], and accounting for the fact that Fontaine writes
$\dG$ multiplicatively. If we write
$\phi^{0}_{\Gscr'}(f)=\Omega_{f}\omega$, then
$$
\nu\phi(u,f)=\Omega_{f}u^{*}\omega
$$
while
$$
\xi\theta(u,f)=f(u)^{*}dZ/(1+Z)\ .
$$
Comparing these formula with the explicit isomorphisms $\xi_{L,\Gscr}$
and $\xi_{L}$ defined in section 1 of [Fon], we see that the
commutativity of the square means that
$$
\Omega_{f}\xi_{L,\Gscr}(u\otimes\omega)=\xi_{L}(f(u)\otimes dZ/(1+Z))\ .
$$
This fact, when combined with   Thm. 1 and Cor. 1  of [Fon], tells us
that
$$
\Omega_{f}\afr_{L}=\afr_{L,\Gscr}\ .
$$
We conclude that
$$
\matrix{
|\Omega_{f}|=p^{-r} & \hbox{with} &
r=\displaystyle\frac{1}{p-1}-\frac{1}{e(q-1)}+\ord_{p}({\cal
D}_{L/\dQ})}
$$
where ${\cal D}_{L/\dQ}$ is the different of the extension $L/\dQ\,$.

To complete the calculation let $t'_{0}$ be our chosen generator for
the $o$-module $T'$.  Some elementary linear algebra using properties
of the different shows that there is a generator $x$ of ${\cal
D}_{L/\dQ}$ such that we have
$$
xt'_{0}=f+f'\ \ \ \hbox{\ in $o_{\dC}\otimes_{\dZ}T'$}
$$
with $f'$ vanishing on the eigenspace in $\dC\otimes_o T$
corresponding to the given embedding $L \subseteq \dC$. This means
that $\phi^{0}_{\Gscr'}(f')=0$ so that
$$
x\Omega_{t'_{0}}=\Omega_{f}\ .
$$
In other words,  the valuation of $\Omega=\Omega_{t'_{0}}$ is
$\ord_{p}(\Omega_{f})-\ord_{p}(x)$ as claimed.

\vfill\eject


{\bf References}

\parindent=23truept

\ref{[Am1]} Amice, Y.: Interpolation $p$-adique. Bull. Soc. math.
France 92, 117-180 (1964)

\ref{[Am2]} Amice, Y.: Duals. Proc. Conf. on $p$-Adic Analysis,
Nijmegen 1978, pp. 1-15

\ref{[BGR]} Bosch, S., G\"untzer, U., Remmert, R: Non-Archimedean Analysis.
Ber-lin-Heidelberg-New York: Springer 1984

\ref{[B-CA]} Bourbaki, N.: Commutative Algebra. Paris: Hermann
1972

\ref{[B-GAL]} Bourbaki, N.: Groupes et alg\`ebres de Lie, Chap. 1-3. Paris:
Hermann 1971, 1972

\ref{[B-VAR]} Bourbaki, N.: Vari\'et\'es diff\'erentielles et
analytiques. Fascicule de r\'e\-sul\-tats. Paris: Hermann 1967

\ref{[Box]} Boxall, J.: $p$-adic interpolation of logarithmic derivatives associated to
certain Lubin-Tate formal groups.  Ann. Inst. Fourier 36 (3), 1-27
(1986)

\ref{[Col]} Colmez, P.: Th\'{e}orie d'Iwasawa des repr\'{e}sentations de de Rham
d'un corps local. Ann. of Math. 148, 485-571 (1998)

\ref{[GKPS]} De Grande-De Kimpe, N., Kakol, J., Perez-Garcia, C.,
Schikhof, W: $p$-adic locally convex inductive limits. In $p$-adic
functional analysis, Proc. Int. Conf. Nijmegen 1996 (Eds.
Schikhof, Perez-Garcia, Kakol), Lect. Notes Pure Appl. Math., vol.
192, pp.159-222. New York: M. Dekker 1997

\ref{[dS]} deShalit, E.: Iwasawa Theory of Elliptic Curves with Complex
Multiplication. Perspectives in Math., vol. 3, Boston: Academic
Press 1987.

\ref{[Fe1]} F\'eaux de Lacroix, C. T.: $p$-adische Distributionen.
Diplomarbeit, K\"oln 1992

\ref{[Fe2]} F\'eaux de Lacroix, C. T.: Einige Resultate \"uber die topologischen
Darstellungen $p$-adischer Liegruppen auf unendlich dimensionalen
Vektorr\"aumen \"uber einem $p$-adischen K\"orper. Thesis, K\"oln
1997, Schriftenreihe Math. Inst. Univ. M\"unster, 3. Serie, Heft
23, pp. 1-111 (1999)

\ref{[Fon]} Fontaine, J.-M.: Formes diff\'{e}rentielles et modules de Tate des
vari\'{e}t\'{e}s ab\'{e}liennes sur les corps locaux.  Invent. math.
65, 379-409 (1982)

\ref{[GR]} Grauert, H., Remmert, R.: Theorie der Steinschen
R\"aume. Berlin-Heidelberg-New York: Springer 1977

\ref{[Kat]} Katz, N.: Formal Groups and $p$-adic Interpolation. In
Ast\'erisque 41-42, 55-65 (1977)

\ref{[Kie]} Kiehl, R.: Theorem A und Theorem B in der
nichtarchimedischen Funktionentheorie. Invent. math. 2, 256-273 (1967)

\ref{[Lan]} Lang, S.: Cyclotomic Fields. Berlin-Heidelberg-New
York: Springer 1978

\ref{[Laz]} Lazard, M.: Les z\'ero des fonctions analytiques d'une
variable sur un corps valu\'e complet. Publ. Math. IHES 14, 47-75
(1962)

\ref{[LT]} Lubin, J., Tate, J.: Formal Complex Multiplication in
Local Fields. Ann. Math. 81, 380-387 (1965)

\ref{[Sch]} Schikhof, W.: Ultrametric calculus. Cambridge Univ.
Press 1984

\ref{[Sh]} Schikhof, W.: Locally convex spaces over nonspherically
complete valued fields I. Bull. Soc. Math. Belg. 38, 187-207
(1986)

\ref{[Sc]} Schneider, P.: $p$-adic representation theory. The 1999
Britton Lectures at McMaster University. Available at
www.uni-muenster.de/math/u/ schneider

\ref{[ST]} Schneider, P., Teitelbaum, J.: Locally analytic
distributions and $p$-adic representation theory, with
applications to $GL_2$. Preprint 1999

\ref{[SI]}  Shiratani, K., Imada, T.:  The exponential
series of the Lubin-Tate groups and $p$-adic interpolation. Mem.
Fac. Sci. Kyushu Univ. Ser A 46 (2), 251-365 (1992)

\ref{[Tat]} Tate, J.: $p$-Divisible Groups. Proc Conf. Local
Fields, Driebergen 1966 (ed. T. Springer), pp. 158-183.
Berlin-Heidelberg-New York: Springer 1967

\bigskip

\parindent=0pt

Peter Schneider\hfill\break Mathematisches Institut\hfill\break
Westf\"alische Wilhelms-Universit\"at M\"unster\hfill\break
Einsteinstr. 62\hfill\break D-48149 M\"unster, Germany\hfill\break
pschnei@math.uni-muenster.de\hfill\break
http://www.uni-muenster.de/math/u/schneider\hfill

\noindent
Jeremy Teitelbaum\hfill\break Department of Mathematics, Statistics,
and Computer Science (M/C 249)\hfill\break University of Illinois at
Chicago\hfill\break 851 S. Morgan St.\hfill\break Chicago, IL 60607,
USA\hfill\break jeremy@uic.edu\hfill\break
http://raphael.math.uic.edu/$\sim$jeremy\hfill

\end